\documentclass[12pt,leqno]{amsart}
\usepackage{amssymb}
\usepackage{amsthm}
\usepackage{amsmath}
\usepackage{amscd}
\usepackage[mathscr]{eucal}
\usepackage[all]{xy}

\numberwithin{equation}{section}

\theoremstyle{plain}
\newtheorem{theorem}{Theorem}[section]

\newtheorem{lemma}[theorem]{Lemma}
\newtheorem{proposition}[theorem]{Proposition}

\theoremstyle{definition}
\newtheorem{definition}[theorem]{Definition}
\newtheorem{remark}[theorem]{Remark}
\newtheorem{example}[theorem]{Example}

\theoremstyle{remark}

\newcommand{\A}{\mathbb{A}}
\newcommand{\E}{\mathbb{E}}
\newcommand{\LL}{\mathbb{L}}
\newcommand{\R}{\mathbb{R}}
\newcommand{\Q}{\mathbb{Q}}

\newcommand{\C}{\mathbb{C}}
\newcommand{\h}{\mathbb{H}}

\newcommand{\W}{\mathbb{W}}

\def\C{{\Bbb C}}
\def\R{{\Bbb R}}

\def\Q{{\Bbb Q}}

\def\eps{\epsilon}

\def\la{\lambda}

\def\8{\infty}
\def\<{\langle}
\def\>{\rangle}

\def\BE{\begin{equation}}
\def\EE{\end{equation}}

\def\8{\infty}

\def\3{\ss}

\newcommand{\G}{\Gamma}
\newcommand{\g}{\gamma}

\newcommand{\x}{\mathbf{x}}

\newcommand{\back}{\backslash}

\newcommand{\calA}{\mathcal{A}}

\newcommand{\calD}{\mathcal{D}}

\newcommand{\calF}{\mathcal{F}}

\newcommand{\calL}{\mathcal{L}}

\newcommand{\calP}{\mathcal{P}}

\newcommand{\calS}{\mathcal{S}}
\newcommand{\calT}{\mathcal{T}}
\newcommand{\calV}{\mathcal{V}}

\newcommand{\tr}{\operatorname{tr}}

\newcommand{\sgn}{\operatorname{sgn}}

\newcommand{\Hom}{\operatorname{Hom}}

\newcommand{\rank}{\operatorname{rank}}

\begin{document}

\title[Cycles with coefficients and modular forms]
{Cycles with local coefficients for orthogonal groups and vector-valued Siegel 
modular forms}

\author[Jens Funke and John Millson]{Jens Funke* and John Millson**}
\thanks{* Partially supported by NSF grant  DMS-0305448 and NSF-grant DMS-0211133 at the Fields Institute, Toronto}
\thanks{** Partially supported by NSF grant DMS-0104006}
\date{May 31, 2004}
\address{Department of Mathematical Sciences, New Mexico State University, 
P.O.~Box 30001, 3MB, Las Cruces, NM 88001, USA}
\email{jfunke@nmsu.edu}
\address{Department of Mathematics, University of Maryland, College Park, MD 
20742, USA}
\email{jjm@math.umd.edu}

\maketitle

\section{Introduction}

The purpose of this paper is to generalize the relation \cite{KM90} between 
intersection numbers of cycles in locally symmetric spaces of orthogonal type 
and Fourier coefficients of Siegel modular forms to the case where the cycles 
have local coefficients. Now the correspondence will involve vector-valued 
Siegel modular forms. 

Let $\underline{V}$ be a non-degenerate quadratic space of dimension $m$ and 
signature $(p,q)$ over $\Q$, for simplicity. The general case of a totally real 
number field is treated in the main body of the paper. We write $V=\underline{V}
(\R)$ for the real points of $\underline{V}$ and let $G = SO_0(V)$. 
Let $G^{\prime}$ denote the nontrivial $2$-fold covering group of the 
symplectic group
$Sp(n,\mathbb{R})$ (the metaplectic group) and $K'$ be the $2$-fold covering
inherited by $U(n)$.  Let $D=G/K$ resp. $D^{\prime} = G'/K'$
be the symmetric space of $G$ resp. $G^{\prime}$. Note that $D' = \h_n$, the 
Siegel upper half space. In what follows
we will choose appropriate (related) arithmetic subgroups  $\G \subset G$ and 
$\G^{\prime}
\subset G^{\prime}$. We let $M = \Gamma\back D$ and $M^{\prime} = \Gamma^
{\prime}
\back D^{\prime}$ be the associated locally symmetric spaces. 
If $M$ is not compact, we let $\overline{M}$ denote the Borel-Serre 
compactification and $\partial M$
denote the Borel-Serre boundary of $\overline{M}$.

We let $\E_n$ denote the holomorphic vector bundle
over $\mathbb{H}_n$ associated to the standard representation
of $U(n)$, i.e., $\E_n = Sp(n,\R) \times_{U(n)} \C^n$. For each dominant 
weight $\la'$ of $U(n)$,
we have the corresponding irreducible representation space 
$S_{\la'}(\mathbb{C}^n)$ of $U(n)$  and the 
associated holomorphic vector bundle $S_{\la'}\E_n$ over $M'$
(see \S \ref{Schur}, for the meaning of the Schur functor $S_{\la'}(\cdot))$. For 
each 
half integer $k/2$ we have a character $det^{k/2}$ of $K'$. Let 
$\mathbb{L}_{k/2}$ be the associated $G'$-homogeneous line bundle over the
Siegel space. For each dominant weight $\lambda$ of $G$, 
we have the corresponding irreducible representation
$S_{[\lambda]}(V)$ of $G$ with highest weight $\lambda$ and the flat
vector bundle $S_{[\lambda]}(\mathcal{V})$ over $M$ with typical fiber
$S_{[\lambda]}(V)$ (see \S \ref{Schur}, for the meaning of the harmonic Schur functor $S_{[\lambda]}(\cdot))$.
 
Let $\lambda$ be a dominant weight for $G$. Let $i(\lambda)$ be the
number of nonzero entries in $\lambda$ when $\lambda$ is expressed
in the coordinates relative to the standard basis $\{\epsilon_i \}$ of \cite
{Bourbaki}, Planche II and IV. Hence we have 
$i(\lambda) \leq [m/2]$.
We will assume (because of the choice of X below in the  construction of our  
cycles $C_X$, see Remark \ref{branching}) that $i(\lambda) \leq p$. Now we 
choose $n$ as in the paragraph
above to be any integer satisfying $i(\lambda)\leq n \leq p$ and choose 
for our highest weight of $U(n)$ corresponding to $\lambda$ the unique dominant 
weight $\lambda^{\prime}$ such that 
$\lambda^{\prime}$ and $\lambda$ {\em have the same nonzero entries},
We note that both weights correspond to the same Young diagram and consequently 
the Schur functors $S_{\lambda^{\prime}}(\cdot)$
and $S_{\lambda}(\cdot)$ are the same so we will not distinguish between them.

The main point of this paper is to use the theta correspondence for  the dual 
pair $(G,G')$ to construct for pair of dominant weights 
$\lambda^{\prime}$ and $\la$ as above  an element
\[
\theta_{nq,[\la]}(\tau,z) \in C^{\infty}(M',S_{\la} \E_n^{\ast}\otimes 
\LL_{-\frac{m}2}) \widehat{\otimes} A^{nq}(M,S_{[\la]} \calV),
\]
($\tau \in \h_n$, $z \in D$) which is closed as a differential form on $M$:
\[
d \theta_{nq,[\la]}(\tau,z) =0.
\]
Our notation is justified since $n$ and $\lambda$ determine $\lambda^{\prime}$.

Hence we obtain an induced element $[\theta_{nq,[\la]}] \in  
C^{\infty}(M',S_{\la} \E^{\ast}_n \otimes \LL_{-\frac{m}2})\otimes  
H^{nq}(M,S_{[\la]} \calV)$. We will say that elements of the above
tensor product are sections of  the holomorphic bundle $S_{\la} \E_n^{\ast}   
\otimes \LL_{-\frac{m}2} $ 
{\em with  coefficients in} $H^{nq}(M, S_{[\la]} \mathcal{V})$. 

Note that the highest weight of the isotropy representation of the homogeneous 
vector 
bundle for the symplectic group coincides (up to a shift) after the addition or
suppression of zeroes with the highest weight of the coefficient system for 
the orthogonal group. This correspondence of the highest weights between 
$O(p,q)$ and $Sp(n,\R)$ agrees with the one obtained by Adams \cite{Ad1}.

On the other hand, we can construct cycles in $M$ as follows. Recall that we 
can realize $D$ as the set of negative $q$-planes in $V$:
\[
D = \{ z \subset V: \dim z =q \quad (\,,\,)|z <0\}.
\]
Then for $\x =(x_1,\dots,x_n) \in \underline{V}^n$ with positive definite inner 
product matrix $(\x,\x) =(x_i,x_j)_{i,j}$, we define a totally geodesic 
submanifolds $D_{\x}$ by
\[
D_{\x} = \{ z \in D: z \perp span(\x) \}.
\]
This gives rise to cycles $C_{\x}$ in $M$ of dimension $(p-n)q$, and by summing 
over all $\x$ in (a coset of) a lattice in $\underline{V}$ such that $\frac{1}
{2}(\x,\x) =\beta>0$, one obtains a composite cycle $C_{\beta}$. For $\beta$ 
positive semidefinite of rank $t \leq n$, there is a similar construction to 
obtain cycles $C_{\beta}$ of dimension $(p-t)q$. We can then assign 
coefficients to these cycles (see \S \ref{cycles} for details) to obtain 
(relative) homology classes\[
C_{\beta,[\la]} \in S_{\la}(\C^n)^{\ast} \otimes H_{(p-t)q}(M, 
\partial M, S_{[\la]} \calV),
\]
i.e., for every vector $w \in S_{\la}(\C^n)$, we obtain a class
\[
C_{\beta,[\la]}(w)  \in H_{(p-t)q}(M, \partial M, S_{[\la]} \calV), 
\]
and for any cohomology class $ \eta \in H_{c}^{(p-n)q}(M,S_{[\la]} \calV)$, the 
natural pairing gives a \emph{vector}
\[
 \langle \eta , C_{\beta,[\la]} \rangle \in  S_{\la}(\C^n)^{\ast}.
\]

In the usual way, we can identify the space of holomorphic sections of the 
bundle $ S_{\la} \E_n^{\ast}\otimes\LL_{-\frac{m}2}$ over (the compactification 
of) $M'$ with $Mod(\G', S_{\la}(\C^n)^{\ast}\otimes \det^{-\frac{m}2})$, the 
space of holomorphic vector-valued
Siegel modular forms for the representation 
$S_{\la}(\C^n)^{\ast} \otimes \det^{-\frac{-m}2} $. 
Here $Mod(\G', S_{\la}(\C^n)^{\ast} \otimes \det^{-\frac{m}2})$ 
is the space of holomorphic functions $f(\tau)$ on $\h_n$ with values in 
$S_{\la}(\C^n)^{\ast} \otimes \det^{-\frac{m}2}$, 
holomorphic at the cusps of $M'$,  such that
\[
f(\g \tau) = (\rho_{\la}^{\ast} \otimes \det {^{-m/2}})( {^t}j(\g,\tau)^{-1}) f
(\tau).
\]
Here $ \rho_{\la}$ is the action of $GL_n(\C)$ on $S_{\la}(\C^n)$ 
and $j(\g,\tau) = c \tau +d$ is the usual automorphy factor for 
$\g = \left( \begin{smallmatrix} a&b\\c&c \end{smallmatrix} \right) \in \G'$. 
Recall that for Siegel modular forms, the Fourier expansion is indexed by 
positive semidefinite $\beta \in Sym_n(\Q)$, and note that the 
$\beta$-th Fourier coefficient of such a form is now a {\em vector} in  
$S_{\la}(\C^n)^{\ast}$.

Our main result is 

\begin{theorem}\label{IMain1}
The cohomology class $[\theta_{nq,[\la]}]$ is a holomorphic Siegel 
modular form for the representation $S_{\la}(\C^n)^{\ast} 
\otimes \det^{-\frac{m}2}$ with coefficients in $  H^{nq}(M,S_{[\la]} \calV)$. 
Moreover, the Fourier expansion of 
$[\theta_{nq,[\la]}](\tau)$ is given by
\[
[\theta_{nq,[\la]}](\tau) =  \sum_{ t=0}^n
\sum_{\substack{\beta \geq 0 \\ \rank \beta = t} } \left( PD(C_{\beta,[\la]})
\cap e_q^{n-t} \right) \, e^{2\pi i tr(\beta \tau)},
\]
where $ PD(C_{\beta,[\la]})$ denotes the Poincar\'e dual class of
$PD(C_{\beta,[\la]})$. Here, for $q$ even, $e_q$ denotes a certain
invariant $q$-form, the Euler form on $D$, and is zero if $q$ is odd.

Furthermore, if $q$ is odd or if $i(\la) =n$, then 
$[\theta_{nq,[\la]}](\tau)$ is a cusp form.

\end{theorem}

This generalizes the main result of \cite{KM90}, where the generating series 
for the special cycles $C_{\beta}$ with trivial coefficients was realized as a 
classical holomorphic Siegel modular form of weight $m/2$.

Pairing with $[\theta_{nq,[\la]}]$ with cohomology and homology defines two 
maps, which we denote both by $\Lambda_{nq,[\la]}$, namely
\[
\Lambda_{nq,[\la]}: H_{c}^{(p-n)q}(M,S_{[\la]} \calV) \longrightarrow 
Mod(\G', S_{\la}^{\ast}\otimes \det{^{-\tfrac{m}2}});
\]
\[
\Lambda_{nq,[\la]}: H_{nq}(M,S_{[\la]} \calV) \longrightarrow 
Mod(\G', S_{\la}^{\ast}\otimes \det{^{-\tfrac{m}2}}).
\]
These pairings give rise to the following two reformulations of Theorem~\ref
{IMain1}:
\begin{theorem}

For any cohomology class $\eta \in H_c^{(p-n)q}(M, S_{[\la]} \calV)$ and 
for any compact cycle $C \in H_{nq}(M,S_{[\la]} \calV)$, the generating series  
\[
 \sum_{t=0}^n \sum_{\beta \geq 0} \langle (\eta \cup  e_q^{n-t}), C_{\beta,
[\la]} \rangle 
 e^{2 \pi i \tr(\beta \tau)}
\]
and 
\[
 \sum_{t=0}^n \sum_{\beta \geq0} \langle C, (C_{\beta,[\la]} \cap e_q^{n-t}) 
\rangle 
 e^{2 \pi i \tr(\beta \tau)}
\] 
define elements in  $Mod(\G', S_{\la}(\C^n)^{\ast}\otimes \det^{-\tfrac{m}2})$.
\end{theorem}

To illustrate our result, we consider the simplest example.

\begin{example}
Consider the weight $\lambda^{\prime}= (\ell,\ell,\cdots,\ell)$ of $U(n)$ 
(so the number or $\ell$'s is $n$). Then $S_{\la}(\C^n) \simeq Sym^{\ell}
(\bigwedge^n (\C^n))$ is one-dimensional, while $S_{[\la]}(V)$ can be realized 
as a summand in the harmonic tensors in $Sym^{\ell} (\bigwedge^n (V)) \subset V^
{\otimes n\ell}$. For $\eta$ a closed rapidly decreasing $S_{[\la]}\calV$-valued 
smooth differential $(p-n)q$-form on $M$, the pairing $\<[\eta], C_{\x,[\la]}
\>$ is given by the period
\[
 \<[\eta], C_{\x,[\la]}\> = \int_{C_{\x}} (\eta,(x_1 \wedge \cdots \wedge x_n)^
{\ell}),
\]
with the bilinear form $(\,,\,)$ on $V$ extended to $V^{\otimes n\ell}$. Then 
the generating series of these periods 
\[
\sum_{\substack{ \x \in L^n \\ (\x,\x)>0}} \left(\int_{C_{\x}} (\eta,(x_1 
\wedge \cdots \wedge x_n)^{\ell}) \right) e^{\pi i tr((\x,\x) \tau)}
\]
is a classical scalar-valued holomorphic Siegel cusp form of weight $ \ell + m/2
$. Here $L$ is (a coset of) an integral lattice in $\underline{V}$. 
\end{example}

For $n=1$, several (sporadic) cases for generating series for periods over 
cycles with nontrivial coefficients as elliptic modular forms were already 
known: For signature $(2,1)$ by Shintani \cite{Shintani}, signature $(2,2)$ 
by Tong \cite{Tong} and Zagier \cite{Zagier}, and for signature $(2,q)$ by 
Oda \cite{Oda} and Rallis and Schiffmann \cite{RS}. For the unitary case of 
$U(p,q)$, see also \cite{TongWang}.

\medskip

We have not tried to prove that the Siegel modular form associated to a
cohomology class $\eta$ or a cycle $C$ is nonzero. However, for the case 
in which $G = SO_0(p,1)$ the nonvanishing of the associated Siegel modular 
form (for a sufficiently deep congruence subgroup depending on $\beta$
and $\lambda$)
follows from \cite{KMCompo} together with \cite{MRaghu}. Indeed first
apply \cite{KMCompo}, Theorem 11.2, to reduce to the case where the
cycle $C_{\beta,
[\lambda]}$ consists of a single component 
$C_{X} \otimes \mathbf{x}_{[f(\lambda)]}$, by passing to
a congruence subgroup, see  
\S \ref{specialcycleswithcoefficients}. Then apply (the proof of)
Theorem 6.4 of \cite{MRaghu} where it is shown that 
for a sufficiently deep congruence subgroup the cycle 
$C_{X} \otimes \mathbf{x}_{[f(\lambda)]}$ is not a boundary.

For general orthogonal groups, the results of J.~S. Li \cite{Li} suggest
that again the Siegel modular form associated to a suitable $\eta$ is nonzero. 
Indeed,
J.~S. Li \cite{Li} has used the theta correspondence (but not our
special kernel $\theta_{nq,[\lambda]}$)  to 
construct non-vanishing cohomology classes for $O(p,q)$ for the above
coefficient systems (with some 
restrictions on $\lambda$). However it is possible that {\em all} \ the above 
cycles
$C_{\beta, [\lambda]}$ are boundaries for some $p,q$ and $\lambda$  This would 
be an 
unexpected development.
\medskip

Finally, we would like to mention our motivation for the present work. We are 
interested in extending the lift $\Lambda_{nq,[\la]}$, say for $\la =0$, the 
trivial coefficient case, to the \emph{full} cohomology $H^{(n-p)q}(M,\C)$. 
This would extend the results of Hirzebruch/Zagier \cite{HZ}, who, for Hilbert 
modular surfaces (essentially $\Q$-rank $1$ for $O(2,2)$), lift the full 
cohomology $H^2(M,\C)$ to obtain generating series for intersection numbers of 
cycles. In this process, cohomology classes and cycles with nontrivial 
coefficients naturally occur, as we now explain.

We study the restriction of $\theta_{nq,0}$ to $\partial(\overline{M})$, which 
is glued together out of faces $e(P)$, one for each $\G$-conjugacy class of 
proper parabolic $\Q$-subgroups $P$ of $G$. We have
\begin{theorem}[\cite{FM1,FM2}]
The theta kernel $\theta_{nq,0}$ extends to $\overline{M}$.  In fact, the 
restriction to $e(P)$ is given by a sum of theta kernels $\theta_{n(q-r),\la}$  
for a nondegenerate subspace $W \subset V$ associated to an orthogonal 
factor of the Levi subgroup of $P$ with values in $S_{\la}(W)$ for certain 
dominant weights $\la$.
\end{theorem}

\medskip

The paper is organized as follows. In \S2, we briefly review homology and 
cohomology with nontrivial coefficients needed for our purposes, while in \S3, 
we review the construction of the finite dimensional representations of $GL_n
(\C)$ and $O(n)$ using the Schur functors $S_{\la}$ and $S_{[\la]}$. We 
introduce the special cycles with coefficients in \S4. In \S5, we give the 
explicit construction of the Schwartz forms $\varphi_{nq,[\la]}$ underlying the 
theta series $\theta_{nq,[\la]}$. We give their fundamental properties and for 
the proofs, we reduce to the case of $n=1$. \S6 is the technical heart of the 
paper, in which we proof the fundamental properties of  $\varphi_{nq,[\la]}$ 
for $n=1$. Our main tool is the Fock model of the Weil representation, which we 
review in the appendix to this paper. Finally, in \S7, we consider the global 
theta series  $\theta_{nq,[\la]}$ and give the proof of the main result.

\medskip
A major part of this work was done while the first named author was a fellow 
at 
the Fields Institute in Toronto in the academic year 02/03. He would like to 
thank the organizers of the special program on Automorphic Forms and the staff 
of the institute for providing such a stimulating environment. We thank Steve 
Kudla for encouraging us to consider the case of an arbitrary dominant weight
in order to produce generating functions for intersection numbers and periods 
that were vector-valued Siegel modular forms.

\section{Homology and cohomology with local coefficients }\label{section2}

In this section we review the facts we need about homology and cohomology
of manifolds (possibly with boundary) with coefficients in a flat bundle 
(``local coefficients'') and ``decomposable cycles''.
We refer the reader to \cite{Hatcher}, page 330 - 336 
for more details.

\subsection{The definition of the groups}
We now define the homology and cohomology groups of $X$ with
coefficients in E. We will do this assuming that 
 $X$ is the underlying space of a connected simplicial complex $K$. 
 We will define the {\em simplicial} homology and cohomology groups 
with values in $E$. By the usual subdivision argument one can prove
that the resulting groups are independent of the triangulation $K$. 

We define a $p$-chain with values in 
$E$ to be a formal sum $\Sigma_{i-1}^m  \sigma_i \otimes s_i$ 
where $\sigma_i$ is an 
{\em oriented} 
$p$-simplex and $s_i$ is a flat section over $\sigma_i$. We denote the group of 
such chains by 
$C_p(X,E)$. Before defining the boundary and coboundary operators we
note that if $t$ is a 
flat section of $E$ over a face
$\tau$ of a simplex $\sigma$ then it extends to a unique flat section
$e_{\sigma,\tau}(t)$ over $\sigma$. Similarly if we have
a flat section $s$ over $\sigma$ it restricts to a flat section 
$r_{\tau,\sigma}(s)$ over $\tau$. Finally if $\sigma = (v_0,\cdots,v_p)$
we define the $i$-th face $\sigma_i$ by 
$\sigma_i = (v_0,\cdots,\hat{v_i},\cdots,v_p)$. Here $\hat{v_i}$ means
the $i$-th vertex has been omitted.

We define the boundary operator $\partial_p : 
C_p(X,E) \longrightarrow C_{p-1}(X,E)$ for $\sigma$ a $p$-simplex
and $s$ a flat section over $X$ by
$$\partial_p (\sigma \otimes s) = 
\sum _{i=0}^p (-1)^i \sigma_i \otimes r_{\sigma_i,\sigma}(s)$$
Then $\partial_{p-1}\circ \partial_p = 0$ and we define the homology groups 
$H_{*}(X,E)$ of $X$ with coefficients
in $E$ in the usual way. These groups depend only on the topological space 
$X$  and the flat bundle $E$.

In a similar way simplicial cohomology groups of $X$ with coefficients in  
$E$ are defined.
A $E$-valued $p$-cochain on $X$ with values in $E$ is a function 
$\alpha$ which assigns to each
$p$-simplex $\sigma $ a flat section of $E$ over $\sigma$.
The coboundary $\delta_p \alpha$ of a $p$-cochain $\alpha$ is defined on a 
$(p+1)$-cochain $\sigma$ by :
$$
\delta_p \alpha (\sigma) = 
\sum _{i=0}^p (-1)^i e_{\sigma,\sigma_i} (\alpha(\sigma_i)).
$$
Then $\delta_{p+1}\circ \delta_p = 0$ and we define the cohomology groups 
$H^{*}(X,E)$
of $X$ with coefficients in $E$ in the usual way.
 
If $A$ is a subspace of $X$, then the complex of simplicial chains with
coefficients in $E|A$ is a subcomplex, and we define the 
relative homology groups $H_.(X,A,E)$ with coefficients in $E$ to be the 
homology
groups of the quotient complex. Similarly, we define the subcomplex of relative
(to $A$) simplicial cochains with coefficients in $E$ to be the complex
of simplicial cochains that vanish on the simplices in $A$ and define
the relative cohomology groups $H^{\bullet}(X,A,E)$ to be the cohomology groups
of the relative cochain complex.

\subsection{Bilinear pairings}

We first define the Kronecker pairing between homology and cohomology with 
local vector bundle
coefficients. Let $E,F$and $G$ be flat bundles over $X$.  Assume that 
$\nu: E \otimes F \longrightarrow G$ is a parallel section of $Hom( E\otimes 
F,G)$.
Let $\alpha$ be a $p$-cochain with coefficients in $E$ and $\sigma \otimes s$ 
be 
a $p$-simplex with coefficients in $F$. 
Then the Kronecker index $<\alpha,\sigma \otimes s >$ is the element of 
$H_0(X,G)$ defined by:
$$< \alpha,\sigma \otimes s > = \nu(\alpha(\sigma) \otimes s).$$
The reader will
verify that the Kronecker index descends to give a bilinear pairing
$$
<\ ,\ >:H^p(X,E) \otimes H_p(X,F) \longrightarrow H_0(X,G).
$$
We note that if $G$ is trivial then $H_0(X,G)\cong G_{x_0}$. In particular,
we get a pairing 
$$
<\ ,\ >:H^p(X,E^*) \otimes H_p(X,E) \longrightarrow \mathbb{R},
$$
which is easily seen to be perfect.
The coefficient pairing $E \otimes F \to G$ also induces cup  products
with local coefficients
$$
\cup :H^p(X,E) \otimes H^q(X,F) \longrightarrow H^{p+q}(X,G)
$$
and cap products with local coefficients (here we assume $m \geq p$)
$$
\cap :H^p(X,E) \otimes H_m (X,F) \longrightarrow H_{m-p}(X,G).
$$
These are defined in the usual way using  the ``front-face'' and ``back-face''
of an ordered simplex and pairing the local coefficients using $\nu$.

\begin{remark} \label{capproduct}
We define the cap product $\alpha \cap \sigma$ for $\alpha$ a $p$-cochain
and $\sigma$ a simplex by making $\alpha$ operate on the {\it back} $p$
face of $\sigma$. 
With this definition the adjoint formula
\begin{equation} \label{adjointformula}
\< \alpha \cup \beta, \sigma \> = \< \alpha ,\beta \cap \sigma \>
\end{equation} 
holds 
(rather than $\< \alpha \cup \beta, \sigma \> = 
\< \beta ,\alpha \cap \sigma \>)$. 
\end{remark} 

The above pairings relativize in a fashion identical to the case
of trivial coefficients. 

Since the proof of Poincar\'e (Lefschetz) duality is a patching
argument of local dualities (see \cite{Hatcher}, p. 245-254), it
goes through for local coefficients as well, and we have

\begin{theorem}
Let $X$ be a compact  oriented manifold with (possibly empty)
boundary and (relative) fundamental class $[X, \partial X]$.
Then we have an isomorphism
$$ 
\mathcal{D}:H^p(X,E)  \longrightarrow H_{n-p}(X,\partial X, E)
$$
given by
$$
\mathcal{D}(\alpha) = \alpha \cap [X,\partial X].
$$
\end{theorem}

\begin{definition}
Suppose $[a] \in H_p(X,\partial X,E)$. We will define the {\em Poincar\'e dual} 
of
$[a]$ to be denoted $PD([a])$ by
$$PD([a]) = \mathcal{D}^{-1}([a]).$$
\end{definition}

We can now define the intersection number of cycles with local coefficients.

\begin{definition}\label{intersectionproduct}
Let $E,F,G$ and $\nu$ be as above and $[a]$ and $[b]$ be homology classes with 
coefficients in $E$ and $F$ respectively. Then we define the intersection
class $[a]\cdot [b] \in H_.(X,G)$ by the formula
$$[a]\cdot [b] = \mathcal{D}( PD([a]) \cup PD([b])).$$
\end{definition}

on our convention in Remark \ref{capproduct} we digress to prove
\begin{lemma}
Using Definition \ref{intersectionproduct} and the convention in Remark \ref
{capproduct} we have  
$$[a]\cdot [b] = \< PD([a]), [b] \>.$$
\end{lemma}

\subsection{Decomposable cycles}

There is a particularly simple construction of cycles with coefficients in $E$.
Let $Y$ be a  compact oriented submanifold with (possibly empty)
boundary $\partial Y \subset \partial X$ of $X$ of codimension $p$ and let
$s$ be a parallel section of the restriction of $E$ to $Y$. Let $[Y,\partial Y]$
denote the relative
fundamental cycle of $Y$ so $[Y,\partial Y]= \Sigma_i \sigma_i$, a sum of 
oriented simplices.

\begin{definition}
$Y\otimes s$ denotes the $(n-p)$-chain with values in $E$ given by
 $$Y = \Sigma_i \sigma_i \otimes s_i$$
where $s_i$ is the value of $s$ on the first vertex of $\sigma_i$.
\end{definition}

\begin{lemma}
$Y\otimes s$ is a relative $n-p$ cycle with coefficients in $E$, called a 
decomposable cycle.
\end{lemma}

For motivation of the term {\it decomposable cycle} we refer the reader to 
\cite{MRaghu}.

\subsection{The de Rham theory of cohomology with local coefficients
and the dual of a decomposable cycle}

In this subsection we recall the de Rham representations of the cohomology
groups $H^*(X,E)$ and of the Poincar\'e dual class $PD(Y \otimes s)$.

>From now on, $X$ will always be smooth manifold.

A differential $p$-form $\omega$ with values in a vector bundle $E$
is a section of the bundle $\bigwedge^p T^*(X) \otimes E$ over $X$. Thus 
$\omega$ assigns to a $p$-tuple of tangent vectors at $x \in X$ a point
in the fiber of $E$ over $x$. 
Suppose now that $E$ admits a flat connection $\nabla$. 
We can then make the graded vector space of smooth $E$-differential forms $A^*
(X,E)$ into a complex by defining
\begin{multline*}
d_{\nabla}(\omega)(X_1,X_2,\cdots,X_{p+1}) = \sum_{i=1}^p (-1)^{i-1} 
\nabla_{X_i}
(\omega(X_1,\cdots,\widehat{X_i},\cdots,X_{p+1})) \\
+ \sum_{i < j} (-1)^{i+j}
\omega([X_i,X_j],X_1,\cdots,\widehat{X_i},\cdots,\widehat{X_j},\cdots,X_
{p+1}) .\end{multline*}
Here $X_i,1\leq i \leq p+1$, is a smooth vector field on $X$.

We now construct a map $\iota$ from $A^p(X,E)$ to the group of
simplicial cochains $C^p(X,E)$ as follows. Let $\omega \in A^p(X,E)$ and
$\sigma$ be a $p$-simplex of $K$. Then in a neighborhood $U$ of $\sigma$ we may
write $\omega = \sum_i \omega_i \otimes s_i$ where the $s_i$'s  are parallel 
sections
of $E|U$ and the $\omega_i$'s are scalar forms. We then define
$$<\iota(\omega),\sigma> = \sum_i (\int_\sigma \omega_i) s_i(v_0).$$

The standard double-complex proof of de Rham's theorem due to Weil, 
see \cite{BottTu}, p. 138, yields

\begin{theorem}
The integration map $\iota:H^*_{de Rham}(X,E) \longrightarrow H^*(X,E)$
is an isomorphism.
\end{theorem}

Finally, we will need that the cohomology class $PD(Y\otimes s)$ has the 
following
representation in de Rham cohomology with coefficients in $E$.

Let $U$ be a tubular neighborhood of the oriented submanifold 
with boundary \ $Y$. We assume (by choosing a Riemannian metric)
that we have a disk bundle $\pi:U \to Y$. Then a Thom form for $Y$
is a closed form $\omega_Y$ where $\omega_Y$ is compactly supported along
the fibers of $\pi$ and has integral one along one and hence all fibers
of $\pi$. It is standard that the extension of $\omega_Y$ to
$X$ by making it zero outside of $U$ represents the Poincar\'e dual of
the class of $[Y,\partial Y]$.  The parallel section $s$ of $E|Y$ extends to a 
parallel section of $E|U$ again denoted $s$. We extend $\omega_Y \otimes s$ to 
$X$ by making it zero outside of $U$. We continue to use the notation
$\omega_Y \otimes s$ for this extended form. We will see below that 
$\omega_Y \otimes s$ represents the Poincar\'e dual of $Y \otimes s$.

If $[a] \in H_p(X,\partial X,E)$, then the de Rham cohomology class $PD([a])$
is the class of $(n-p)$--forms characterized by the property that 
if $a$ is a simplicial cycle representing $[a]$, then for any $E^*$ valued 
$p$-form $\eta$ vanishing on $\partial X$ we have 

\begin{equation*}
\int_X \eta \wedge PD([a]) = \int_a \eta.
\end{equation*}

\begin{remark}
In abstract terms the above equation is 
$$\<[\eta] \cup PD([a]), [X,\partial X] \> = \<[\eta],[a]\>.$$
Since the expression on the right-hand side of this formula
is equal to $\<[\eta], PD([a]) \cap [X,\partial X] \>$ our definition
of Poincar\'e dual amounts to assuming the adjoint formula, \eqref
{adjointformula}, and hence that assuming the ``back
p face '' definition of the cap product.
\end{remark}
 
\begin{lemma}\label{Poincaredual}
Then the de Rham cohomology class Poincar\'e dual to the 
cycle with coefficients $Y \otimes s$
is represented by the bundle--valued form $\omega_Y \otimes s$.
\end{lemma}
\proof
We need to prove that for any $E^*$-valued closed $(n-p)$-form $\eta$ 
vanishing on $\partial X$ we 
have
$$\int_X \eta \wedge \omega_Y \otimes s = \int_{Y \otimes s} \eta = 
\int_Y \<\eta,s\>.$$
But 
$$\int_X \eta \wedge \omega_Y \otimes s = \int_X \<\eta,s\>\wedge \omega.$$
But since $s$ is parallel on $U$ the scalar form $\<\eta,s\>$ is closed,
and the lemma follows because $\omega_Y$ is the Poincar\'e dual to 
$[Y, \partial Y]$.
\qed

\section{Finite dimensional representations of $GL(n)$ and $O(n)$}\label{Schur}

In this section we will review the construction of the irreducible
finite dimensional (polynomial) representations of $GL(U)$ (resp.
$O(U)$) where $U$ is a complex vector space of dimension $n$
(resp. a finite dimensional complex vector space of dimension $n$ equipped 
with a non-degenerate symmetric bilinear form $( \ ,\ )$). 

\subsection{Representations of the general linear group}

\subsubsection{Schur functors}

We recall that the symmetric group $S_{\ell}$ acts on the $\ell$--fold 
tensor product $T^{\ell}(U)$ according to the rule
that $s \in S_{\ell}$ acts on a decomposable $v_1 \otimes \cdots v_{\ell}$
by moving $v_i$ to the $s(i)$--th position.
Let $\lambda = (b_1,b_2,\cdots, b_n)$ be a partition of $\ell$. We assume
that the $b_i$'s are arranged in decreasing order. We will use $D(\lambda)$
to denote the Young diagram associated to $\lambda$.

For more details on what follows, see \cite{FultonHarris}, \S 4.2 and \S 6.1,
\cite{GoodmanWallach}, \S 9.3.1--9.3.4 and \cite{Boerner}, Ch.~V, \S 5.

\paragraph{\bf Standard fillings and the associated projections}

\smallskip

\begin{definition}
A standard filling $t(\lambda)$ of the Young diagram 
$D(\lambda)$  by the elements of the set 
$[\ell]=\{1,2, \cdots,  \ell \}$ 
is an assignment of each of the numbers in $[\ell]$ to a box of $D(\lambda)$ so
that the entries in each row strictly increase when read from left to right 
and the entries in each column strictly increase when read from top to bottom.
We will denote the set of standard fillings of $D(\lambda)$ by $S(\lambda)$.
A Young diagram equipped with a standard filling will be called a standard
tableau. 
\end{definition}

We let $t_0(\lambda)$ be the standard filling that
assigns $1,2,\cdots,\ell$  from left to right starting with the first row
then moving to the second row etc.

We now recall the projection in $End(T^{\ell}(U))$ associated to
a standard tableau $T$ with $\ell$ boxes corresponding to a standard filling
$t(\lambda)$ of a Young diagram $D(\lambda)$. 
Let $P$ (resp. $Q$) be the group preserving the rows (resp. 
columns) of $T$. Define elements of the group ring of $S_{\ell}$
by  $r_{t(\lambda)} = c_1 \sum_P p$ and $c_{t(\lambda)} = c_2 \sum_Q \epsilon
(q) q$
where $c_1 = 1/|P|$ and $c_2 = 1/|Q|$, so $r_{t(\lambda)}$ and $c_{t(\lambda)}$
are idempotents. We let $\mathcal{P}$ (resp. $\mathcal{Q}$) be the projections
operating on $T^{\ell}(U)$ obtained by acting by $r_{t(\lambda)}$ 
(resp. $c_{t(\lambda)}$).  We put $s_{t(\lambda)} = c_3 c_{t(\lambda)} \cdot 
r_{t(\lambda)}$
(product in the group ring) where  $c_3$ is chosen so that 
$s_{t(\lambda)}$ is an idempotent, see \cite{FultonHarris}, Lemma 4.26.

\begin{remark}
We have abused notation by not indicating the dependence of 
$\mathcal{Q}$ and $\mathcal{P}$ on the standard filling $t(\lambda)$.
We will correct both these abuses by letting $\pi_{t(\lambda)}$ denote
the projector obtained by correctly normalizing the previous product.
\end{remark}

We now have

\begin{theorem}

We have a direct sum decomposition
\[
T^{\ell} (U) = \bigoplus_{\la \in \calP(\ell)} \bigoplus_{t(\la) \in 
S(\la)} s_{t(\la)} \left( T^{\ell} (U) \right),
\]
where $\calP(\ell)$ denotes the set of partitions of $\ell$.
\end{theorem}

Furthermore we have, \cite{GoodmanWallach}, Theorem 9.3.9,

\begin{theorem}
For every standard filling $\lambda$, the $GL(V)$--module 
$s_{t(\la)} \left( T^{\ell} (U) \right)$ is irreducible with
highest weight $\lambda$.
\end{theorem}

\begin{remark}
If $t^{\prime}(\lambda)$ is another standard filling then the permutation
relating the two fillings induces an isomorphism of $\pi_{t(\lambda)}(T^{\ell}
(U))$
and $\pi_{t^{\prime}(\lambda)}(T^{\ell}(U))$.
\end{remark}

For concreteness we will define the Schur functor $S_{\lambda}( \cdot)$
 by choosing $t(\lambda) = t_0(\lambda)$, whence
$S_{\lambda}(U) =  s_{t_0(\la)} \left( T^{\ell} (U) \right)$. We
 obtain projections $\pi_{\lambda}$
\[
\pi_{\la}: T^{\ell} (U) \longrightarrow S_{\la} (U) 
\]
and inclusions
\[
\iota_{\la}: S_{\la} (U)  \longrightarrow T^{\ell} (U).\]

\paragraph{\bf Semistandard fillings and the associated basis of 
$S_{\lambda}(U)$}

\smallskip

\begin{definition}
A semistandard filling of $D(\lambda)$ by the set $[n] =\{1,2, \cdots, n\}$ 
is an assignment of the numbers in $[n]$ to the boxes of $D(\lambda)$ such
that the numbers in each row increase weakly and the numbers in each
column increase strictly. We let $SS(\lambda,n)$ denote the set of
semistandard fillings of $D(\lambda)$ by the elements of the set $[n]$.
\end{definition}

We will also need

\begin{definition}
Suppose $\x = (x_1,\cdots,x_n) \in U^n$  and
$f(\lambda) \in SS(\lambda,n)$. Suppose $a_{ij}$ is the $j$-th entry in the 
$i$-th column. Then $x_{f(\lambda)}$, the word in $\x$ corresponding
to $f(\lambda)$, is defined by
$$
\x_{f(\lambda)} = x_{a_{11}} \otimes x_{a_{12}} \otimes \cdots \otimes
x_{a_{kb_k}}.
$$
\end{definition}

We have

\begin{theorem}\label{labasis}
Let $u_1,\cdots,u_n$ be a basis for $U$ and let $u=(u_1,\dots,u_n)$.
Then the set of vectors $\{ \pi_{t(\lambda)}(u_{f(\lambda)}): f\in SS
(\lambda,n) \}$ is a basis
for $\pi_{t(\lambda)}(T^d(U))$.
\end{theorem}

For a simple proof of this theorem see \cite{Boerner}, Theorem 5.3. Boerner
proves the theorem using the idempotent $ c_3 \mathcal{P} \mathcal{Q}$
on $T^{\ell}(U)$ (actually, he considers $c_3  \mathcal{Q}\mathcal{P}$ on
$T^{\ell}(U^{\ast})$), but his proof can be easily modified to give the
theorem above.

\subsection{Representations of the orthogonal group}

\bigskip

\subsubsection{The harmonic Schur functors}

We will follow \cite{FultonHarris} in our description of the harmonic
Schur functor $U \to S_{[\lambda]}(U)$ on an $n$-dimensional non degenerate 
quadratic space $(U, (\ ,\ ))$ corresponding to a partition $\lambda$.

\smallskip

\paragraph{\bf The harmonic projection}

\smallskip

Suppose that $V_{\lambda}$ is the irreducible representation of $SO(U)$ 
with highest weight $\lambda = (b_1,b_2,\cdots,b_k)$ where $k= [\frac{n}{2}]$.
We will abuse notation and use $\lambda$ to denote the corresponding partition
of $\ell = \sum b_i$. We extend the quadratic form $(\ ,\ )$  to $T^{\ell}(U)$
as the $\ell$-fold tensor product and note that the action of $S_\ell$ on
$T^{\ell}(U)$ is by isometries. For each pair $I = (i,j)$ of integers
between $1$ and $\ell$ we define the contraction operator $C_I
:\otimes^\ell U \to \otimes^{\ell-2} U$ by
$$C_I(v_1\otimes \cdots v_{\ell}) = \sum_k (v_i,e_k)(v_j,e_k) v_1\otimes \cdots
\otimes \widehat{v_i} \otimes \cdots \otimes \widehat{v_j}\otimes \cdots 
\otimes v_{\ell}$$
where $\{e_1,\cdots,e_n\}$ is an orthonormal basis for $(\ ,\ )$.
We also define the expansion operator 
$A_I:\otimes^{\ell-2} U \to \otimes^\ell U$ to be the adjoint of $C_I$, 
that is the operator that inserts the  (dual of)
the form $(\ ,\ )$ into the $(i,j)$-th spots. 
 We define the harmonic 
$\ell$-tensors, to be
denoted $U^{[\ell]}$, to be the kernel of all the contractions $C_I$.
Following \cite{FultonHarris}, p. 263, we define the subspace 
$U^{[\ell]}_{\ell-2r}$ of $U^{\otimes \ell}$ by

$$U^{[\ell]}_{\ell-2r} = \sum A_{I_1} \circ \cdots 
A_{I_r}U^{[\ell -2r]}.$$

Carrying over the proof of \cite{FultonHarris}, Lemma 17.15 (and the
exercise that follows it) from the symplectic case to the orthogonal case
we have

\begin{lemma}
We have a direct sum, orthogonal for $(\ ,\ )$,
$$T^{\ell}(U) = U^{[\ell]} \oplus \oplus_{r=1}^{[\frac{\ell}{2}]} 
U^{[\ell]}_{\ell - 2r}.$$
\end{lemma}

We define the harmonic projection $\mathcal{H}:T^{\ell}(U) \to U^{[\ell]}$ to 
be the 
orthogonal projection onto the harmonic $\ell$-tensors $U^{[\ell]}$. The space 
of
harmonic $\ell$-tensors $U^{[\ell]}$ is invariant under the action of
$S_\ell$. Consequently we may apply the idempotents in the group algebra
of $S_\ell$ corresponding to partitions to further decompose $U^{[\ell]}$
as an $O(U)$--module.

\smallskip

\smallskip
\paragraph{\bf The harmonic Schur functors}
\smallskip

Again following \cite{FultonHarris}, p. 296, we define the harmonic Schur 
functor $S_{[\lambda]}U$ as follows.

\begin{definition}
$$S_{[\lambda]}(U)= \mathcal{H}\mathcal{Q}\mathcal{P} T^{\ell}(U) = \mathcal{H}
S_{\lambda}(U).$$
\end{definition}

We then have see \cite{FultonHarris}, Theorem 19.22,

\begin{theorem} 
The $O(U)$ module $S_{[\lambda]}(U)$ is irreducible with highest weight
$\lambda$.
\end{theorem}

\begin{definition}
We write $\pi_{[\la]}= \mathcal{H} \circ \pi_{\la}$ for the projection
from $T^{\ell}(U)$ onto  $S_{[\lambda]}(U)$. For a semistandard filling $f(\la)
$, we also set for $\x = (x_1,\dots,x_n)$,
\[
\x_{[f(\la)]} = \pi_{[\la]} \x_{f(\la)} \in S_{[\la]}(U).
\]
\end{definition}

In what follows, we will need the following

\begin{lemma}
\begin{enumerate}
\item[(i)] $\mathcal{H},\mathcal{P}$ and $\mathcal{Q}$ are self-adjoint
relative to $(\ ,\ )$.
\item[(ii)] $\mathcal{H}$ commutes with $\mathcal{P}$ and $\mathcal{Q}$.
\end{enumerate}
\end{lemma}

\proof
It is clear that $\mathcal{H}$ is self-adjoint. The arguments for $\mathcal{P}$ 
and
$\mathcal{Q}$ are the same. We give the one for $\mathcal{Q}$.
We will use the symbol $q$ to denote both the element $q \in Q$
and the corresponding operator on $T^{\ell}(V)$. Since $q$
is an isometry we have $q^* = q^{-1}$. Hence we have
$$
\mathcal{Q}^* = \sum \epsilon(q)q^* = \sum \epsilon(q^{-1}) q^{-1}
= \mathcal{Q}.
$$
To prove that $\mathcal{H}$ commutes with $\mathcal{P}$ and $\mathcal{Q}$
it suffices to prove that $\mathcal{H}$ commutes with every element 
$g \in S_{\ell}$. But $S_{\ell}$ acts by isometries and preserves $U^{[\ell]}$.
Consequently it commutes with orthogonal projection on $U^{[\ell]}$.
\qed

\section{Special cycles with local coefficients}\label{cycles}
 
\subsection{Arithmetic quotients for orthogonal groups}

Let $\mathbb{K}$ be a totally real number field with Archimedean places
$v_1,\dots,v_r$ and associated embeddings $\la_1,\dots,\la_r$ and let be
$\mathcal{O}$ its ring of algebraic
integers. Let $\underline{V}$ be an oriented vector space over 
$\mathbb{K}$ of dimension $m \geq 3$ with a non-degenerate bilinear form 
$(\,,\,)$  and let $V$ be the completion of 
$\underline{V}$ at $v_1$. We assume that the associated quadratic form
has signature $(p,q)$ at the completion $v_1$ and is 
positive definite at all other completions. Finally, we let $L$ be an integral 
lattice in $\underline{V}$ and $L^{\#} \supseteq L$ its dual lattice.

Let $\underline{G}$ be
the algebraic group whose $\mathbb{K}$-points is the group of
orientation preserving isometries of determinant $1$ of the form $(\,,\,)$ and 
let  $G:=\underline{G}(\mathbb{R})$ its real points. We let $\Phi = \underline
{G}(\mathcal{O})$ be the subgroup of $\underline{G}(\mathbb{K})$ consisting of 
those elements that take $L$ into itself. 
We let $\mathfrak{b}$ be an ideal in 
$\mathcal{O}$ and let $\Gamma = \Gamma(\mathfrak{b})$ be the congruence
subgroup of $\Phi$ of level $\mathfrak{b}$ (that is, the elements of
$\Phi$ that are congruent to the identity modulo $\mathfrak{b}$). We fix a 
congruence condition $h \in L^n$ once and for all and note that $\G$ operates 
on the coset $ h + \mathfrak{b}L^n$.

We realize the symmetric space associated to $V$ as the set of negative $q$-
planes in $V$:
\begin{equation*}
D \simeq \{ z \subset  V;\, \dim z =q \quad \text{and} \quad
(\,,\,)|_z < 0\}.
\end{equation*} 
We denote the base point of $D$ by $z_0$, and we have $D \simeq G/K$, where $K$ 
is the maximal compact subgroup of $G$ stabilizing $z_0$. Also note $\dim_{\R} 
D = pq$. For $z \in D$, we write $(\,,\,)_z$ for the associated majorant.
Finally, we write
\[
M =\G \back D
\]
for the locally symmetric space.

\subsection{Special cycles with trivial coefficients}

Let ${\bf x}=\{x_1,x_2,\cdots,x_n\} \in \underline{V}^n$ be an $n$-tuple of 
$\mathbb{K}$--rational vectors. We let $\underline{X}$ be the span
of ${\bf x}$ and let $X$ be the completion of $\underline{X}$ at $\la_1$. 
We write $(\x,\x)$ for the $n$ by $n$ matrix with $ij$--th entry equal to
$(x_i,x_j)$. We call ${\bf x}$ nondegenerate if $\rank (\la_i \x, \la_i \x) = 
\dim \underline{X}$ for all $i$ and nonsingular if  $\rank (\la_i \x, \la_i \x)
=n$.

Assume $\x$ is nondegenerate with $\dim \underline{X} = t \leq n$ such that 
$(\ ,\ )|\underline{X}$ is positive definite. Let $r_X$ be the isometric 
involution of $V$ given by
$$r_X(v) =
\begin{cases}
-v &   \text{if $ v\in X$} \\
v  &  \text{if $v\in X^{\perp}$}.
\end{cases}
$$
We define the totally geodesic subsymmetric space $D_X$ by
$$D_X = \{z \in D : (z,x_i) = 0, 1\leq i \leq n\}.$$
Then $D_X$ is the fixed-point set of $r_X$ acting on $D$ and has codimension 
$(n-t)q$ in $D$. We orient $D_X$ as in \cite{KM90}, p.130-131. We also define 
subgroups $G_X$ (resp. $\Gamma_X$) to be the stabilizer
in $G$ (resp. in $\Gamma$) of the subspace $X$. We define $G^{\prime}_X
\subset G_X$ to be the subgroup that acts trivially on $X$, and put $\G'_X = \G 
\cap G'_X$.

\begin{theorem}\label{embeddedcycles} 
There exists a congruence subgroup $\Gamma := \Gamma(\mathfrak{b})$ of
$\Phi$ such that

\begin{enumerate}
\item $M =\Gamma\backslash D$ is an orientable 
 manifold of dimension $pq$ with finite volume, and 
\item
for all $X$ as above,  the image $C_X$ of $D_X$ in $M$ is the 
quotient $\Gamma_X\backslash D_X$ and defines a properly embedded orientable 
submanifold of 
codimension $(n-t)q$. 
\end{enumerate}
\end{theorem}

The theorem will be a consequence of the existence of a ``neat''
congruence subgroup. We recall the definition of a neat subgroup of $\Gamma$.

\begin{definition}[\cite{Borel},p.~117]\label{neatsubgroup}

 An element $g \in G$ is neat if the subgroup of $\mathbb{C}^*$
generated by the eigenvalues of $g$ is torsion free.
In particular, if a root of unity $z$ is an eigenvalue of a neat element 
then $z=1$. A subgroup $\Gamma \subset G$ is neat if all the elements in 
$\Gamma$
are neat.

\end{definition}

We have 

\begin{proposition} [Proposition~17.4, \cite{Borel}]\label{existenceofneat
subgroup}
Let $G$ be an algebraic group defined over $\mathbb{Q}$ and $\Gamma$
an arithmetic subgroup. Then $\Gamma$ admits a neat congruence subgroup.
\end{proposition}

Theorem \ref{embeddedcycles} is an immediate consequence of the
following

\begin{lemma} \label{trivialaction}
If $\Gamma$ is a neat subgroup, then $\Gamma_X$
acts trivially on $X$, i.e., $\G'_X = \G_X$.
\end{lemma}
\proof
We have a projection map $p_X :\Gamma_X \to O(X_1) \times O(X_2) \times \cdots
\times O(X_r)$. Here by $X_i$ we mean the $i$-th completion of $X$.
The $i$-th completion of $(\ ,\ )$ restricted to $X_i$ is positive
definite for $1\leq i \leq r$. Furthermore the splitting $V=X \oplus
X^{\perp}$ is defined over $\mathbb{K}$. Thus the diagonal embedding
of the intersection $L_X = L \cap X$ is a lattice in $\oplus_{i=0}^r 
X_i$ which is invariant under $p_X(\Gamma_X)$. Hence $p_X(\Gamma_X)$
is a discrete subgroup of a compact group hence a finite group. Hence
if $\gamma \in p_X(\Gamma_X)$, then all eigenvalues of $\gamma$ are roots
of unity. Since $\Gamma$ is neat all eigenvalues must be $1$ and the
lemma follows.
\qed

We will later need

\begin{definition} 
The Riemannian exponential map from the total space of the normal
bundle of $D_X$ to $D$ induces a fiber bundle $\pi_X: D \to D_X$ with
totally geodesic fibers. The map $\pi_X$ induces a quotient fibering
$\pi_X: \G_X \back D \to \G_X \back D_X = C_X$, see \cite{KMCompo}. A
Thom form  $\Phi_X$ for the cycle $C_X$ is a closed integrable $(n-t)q$-form on
$\G_X \back D$ such that the integral of $\Phi_X$ over each fiber of
$\pi_X$ is $1$. In particular, $\Phi_X$ is a Poincar\'e dual form for the cycle
$C_X$ in the noncompact submanifold $\G_X \back D$.
\end{definition}

Occasionally, we will also write $C_{\x}$ ($D_{\x}$) for $C_X$ ($D_X$).

We introduce composite cycles as follows. For $\beta \in Sym_n(\mathbb{K})$, we 
set
\[ 
\Omega_{\beta} = \{ \x \in \underline{V}^n: \tfrac12 (\x,\x) = \beta \}
\]
and 
\[
 \Omega^c_{\beta} = \{ \x \in \Omega_{\beta}: \dim \la_i \underline{X} = \rank 
\beta \quad \text{for all $i$}\}. 
\]
We put 
\[
\mathcal{L}_{\beta} = \mathcal{L}_{\beta}(h,\mathfrak{b}) =  
(h + \mathfrak{b}L^n) \cap \Omega_{\beta}.
\]
Then $\G$ acts on $\mathcal{L}^c_{\beta} = \mathcal{L}_{\beta} \cap \Omega^c_
{\beta}$ with finitely many orbits and for $\beta$ positive semidefinite (i.e., 
$\la_i(\beta) \geq 0$ for all $i$), we define
\[ 
C_{\beta} = \sum_{\x \in \G \back \mathcal{L}^c_{\beta}} C_X. 
\]

\subsection{Special cycles with nontrivial coefficients}
\label{specialcycleswithcoefficients}

We now want to promote $C_X$ to a (decomposable) cycle with coefficients
for appropriate coefficient systems $W$ by finding a nonzero 
parallel section of $\mathcal{W}|C_X$. Note that it is
enough to find any $\G_X$-fixed vector $w\in W$ since such a vector $w$ gives 
rise to a parallel section $s_w$  of $\mathcal{W}|C_X$ in the usual
way. Namely, for $z \in C_X$, the section $s_w$ for the bundle  $C_X \times_
{\G_X} W \to C_X$ is given by $s_w(z) = (z,w)$. Thus $s_w$ is constant, hence
parallel. Furthermore, for such a vector $w$, we write $C_X \otimes w$ for
$C_X \otimes s_w$.

The key point for us in constructing parallel sections is Lemma
\ref{trivialaction}. Namely, the components $x_1,\cdots,x_n$ of ${\bf
  x}$ are all fixed by $\Gamma_X = \G_X'$, hence any tensor word in
these components will be fixed by $\G_{X}$.

\begin{definition}
For $f(\la)$ a semistandard filling for $D(\la)$, we define special cycles with
coefficients in $S_{[\la]}(V)$ by setting
\begin{gather*}
C_{\x,[f(\la)]}= C_X \otimes (\la_1\x)_{[f(\la)]}.
\end{gather*}
We also define composite cycles 
$C_{\beta,[f(\la)]}$ analogously as before.
\end{definition}

To lighten the notation, we will write in the following $\x_{f(\la)}$ and $\x_
{[f{\la}]}$ for $(\la_1\x)_{f(\la)}$ and  $(\la_1\x)_{[f(\la)]}$.

\medskip

However, there is an obstruction to the construction of nonzero sections.

\begin{remark}[\cite{MRaghu}, Proposition~4.3] \label{branching}
Let $\la$ be the highest weight of $W$ and let $i(\la)$ be the number
of nonzero entries in $\la$ (so $i(\la)$ is the number of rows
in the associated partition). Then
\begin{equation*} 
dim(X) \geq i(\la)
\end{equation*}
is a necessary condition for the existence of a $G^{\prime}_X$--invariant
vector in $W$ i.e. to finding a nonzero parallel section of the restriction
of the flat vector bundle $\mathcal{W}$ to the cycle $C_X$.
\end{remark} 

On the other hand, for $  dim(X) \geq i(\la)$, we do have nonzero  parallel 
sections along the submanifold $C_X$.

\begin{theorem}[ \cite{MRaghu},Theorem 4.13]
For a weight $\la=(b_1,\dots,b_{[\frac{m}{2}]})$, assume $i(\la) = k \leq n$. 
 Let $f_0(\la)$ be the semistandard filling that puts $1$'s in the first row of 
$D(\mu)$, $2$'s in the second row etc.. Furthermore assume that for  $\x = 
(x_1,\dots,x_n)$, the first $k$ vectors $x_1,\dots,x_k$ are linearly 
independent and satisfy
$(x_i,x_j) = 0, i \neq j$. Then
\[
\x_{[f_0(\la)]} = \mathcal{H} \pi_{\la}(\x_{f_0(\la)}) = 
\mathcal{H} \mathcal{Q}
 x_1^{\otimes b_1}\otimes \cdots \otimes x_k^{\otimes b_k} 
\]
 is a nonzero $\G_X^{\prime}$--invariant in $S_{[\la]}(V)$.
\end{theorem}


For later use, we record (by an analog of Lemma~\ref{Poincaredual})

\begin{lemma}\label{PDform}
Let $\eta$ be  a rapidly decreasing $\mathcal{
  S_{[\la]}({V})  }$--valued closed $(p-n)q$ form on $M$. 
If $\Phi_X$ denotes a Thom form for the cycle $C_X$, then $\Phi_X
\otimes \x_{[f(\la)]}$ satisfies
\[
\int_M \eta \wedge \left(\Phi_X\otimes \x_{[f(\la)]} \right) =
\int_{C_{\x,[f(\la)]}} \eta = \int_{C_X} (\eta, \x_{[f(\la)]}).
\]

\end{lemma}
\subsection{Cycle-valued homomorphisms on $T^{\ell}(\Q^n)$}

We now construct composite cycles $C_{\beta,[\la]}$ 
which are homomorphisms from $S_{\la}(\mathbb{Q}^n)$ to  $H_{\bullet}(M, S_{[\la]}
(\mathcal{V})$. 

\begin{definition}
We define elements $C_{X,[\la]}(\cdot)$
of $Hom(S_{\la}(\mathbb{Q}^n),H_{\bullet}(M, S_{[\la]}(\mathcal{V}))$ 
by
\begin{gather*} 
C_{\x,[\la]} (\eps_{f(\la)}) = C_X \otimes \x_{[f(\la)]}.
\end{gather*}
(Note the map automatically factors  through $S_{\la}(\mathbb{Q}^n)$).
We then have composite cycles $C_{\beta,[\la]}(\cdot)$ as before by summing 
over all $\x \in \mathcal{L}_{\beta}$.

\end{definition}

Finally note that, if $\eta$ is   a rapidly decreasing ${ S_{[\la]}(\mathcal
{V})  }$--valued closed $(p-n)q$ form
on $M$, then the period $\int_{C_{\x,[\la]}} \eta$ is the linear functional
on $S_{\la}(\mathbb{C}^n)$ given by
\begin{equation}
\left(\int_{C_{\x,[\la]}} \eta\right) (\pi_{\la} \epsilon_{f(\la)})
= \int_{C_X} (\eta, \x_{[f(\la)]}).
\end{equation}

\section{Special Schwartz forms}

In this section, we will explicitly construct the Schwartz form $\varphi_{nq,
[\la]}$ needed to construct the cohomology class 
$[\theta_{nq,[\la]}](\tau,z)$ alluded to in the introduction. As in the 
introduction we will choose a pair of highest weights
$\lambda$ for $G$ and $\lambda^{\prime}$ for $U(n)$ which have the same nonzero 
entries. We let $\ell$ be the sum of
the entries of $\lambda$ (equals the sum of the entries of $\lambda^{\prime}$).

\subsection{A double complex for the Weil representation}

In this section, $V$ will denote a real quadratic space of dimension $m$ and 
signature $(p,q)$. We write $\calS(V^n)$ for the space of (complex-valued) 
Schwartz functions on $V^n$. We denote by $G'= Mp(n,\R)$ the metaplectic cover 
of the symplectic group $Sp(n,\R)$ and let $K'$ be the inverse image of the 
standard maximal compact $U(n) \subset Sp(n,\R)$ under the covering map $Mp
(n,\R) \rightarrow Sp(n,\R)$. Note that $K'$ admits a character $\det^{1/2}$, 
i.e., its square descends to the determinant character of $U(n)$. The embedding 
of $U(n)$ into $Sp(n,\R)$ is given by $A + iB \mapsto 
\left( \begin{smallmatrix}
A & B \\
-B &A
\end{smallmatrix}\right)$.
We let $\omega = \omega_V$ be the Schr\"odinger model of the (restriction of 
the) Weil representation of $G' \times O(V)$ acting on $\calS(V^n)$ associated 
to the additive character $t \mapsto e^{2\pi i t}$. 

We let $\h_n = \{ \tau = u+iv \in Sym_n(\C): \; v >0 \} \simeq Sp(n,\R)/U(n) $ 
be the Siegel upper half space of genus $n$. We write $\mathfrak{g}'$ and 
$\mathfrak{k}'$ for the complexified Lie algebra of $Sp(n,\R)$ and $U(n)$ 
respectively. We write the Cartan decomposition as $\mathfrak{g}' = \mathfrak
{k}' \oplus \mathfrak{p}'$, and write $\mathfrak{p}' = \mathfrak{p}^+ \oplus  
\mathfrak{p}^-$ for the decomposition of the tangent space of the base point 
$i1_n$ into the holomorphic and anti-holomorphic tangent spaces.
We let $\bar{Z}_{j}$, $1 \leq j \leq n(n+1)/2$ be a basis of
$\mathfrak{p}^-$ and let $\bar{\eta}_{j}$ be the dual basis. We let
$\C(\chi_{m/2})$ be the $1$-dimensional representation $\det^{m/2}$ of
$K'$. We write $W_{\ell} = T^{\ell} (\C^n) \otimes \C(\chi_{m/2})$
considered as a representation of $K'$ and let $\mathcal{W}_{\ell}$ be
the $G'$-homogeneous vector bundle over $\h_n$ associated to
$W_{\ell}$. We also define $W_{\lambda^{\prime}}$ and $\mathcal{W}_{\lambda^
{\prime}}$ in the
same way using $S_{\lambda^{\prime}}(\C^n)$ instead.

We pick an orthogonal  basis $\{e_i\}$ of $V$ such that
$(e_{\alpha},e_{\alpha}) = 1$ for $\alpha = 1,\dots,p$ and
$(e_{\mu},e_{\mu}) = -1$ for $\mu = p+1,\dots,p+q$. In what follows, we will 
use ``early'' Greek letters (typically $\alpha$ and $\beta$) as subscripts to 
denote indices between $1$ and $p$ (for the ``positive'' variables) and 
``late'' ones (typically $\mu$ and $\nu$) to denote indices between $p+1$ and 
$p+q$ (for the ``negatives'' ones).

Let $\mathfrak{g}$ be the Lie algebra of $G$ and $\mathfrak{g} =
 \mathfrak{p} + \mathfrak{k}$ its Cartan decomposition, where $Lie(K)
 = \mathfrak{k}$. Then $\mathfrak{p} \simeq \mathfrak{g}/  \mathfrak{k}$ is 
isomorphic to the tangent space at the base point of $D \simeq G/K$.
We denote by $X_{\alpha\mu}$ ($ 1 \leq \alpha \leq p$, $p+1 \leq \mu \leq p+q$) 
the elements of the standard basis of $\mathfrak{p}$ induced by the basis 
$\{e_i\}$ of $V$, i.e.,
\[
X_{\alpha\mu}(e_i) = 
\begin{cases}
e_{\mu}, & \text{if $i =\alpha$} \\
e_{\alpha}, &  \text{if $i = \mu$} \\
0, & \text{otherwise.}
\end{cases}
\]
We let $\omega_{\alpha\mu} \in \mathfrak{p}^{\ast}$ be the elements of the 
associated dual basis. Finally, we let  $\mathcal{A}^k(D)$ be the space of 
(complex-valued) differential $k$ forms on $D$.

\smallskip

We consider the graded associative algebra 
\[
C = \bigoplus_{i,j,\ell \geq 0} C^{i,j}_{\ell},
\]
where
\[
C^{i,j}_{\ell} = \left[ {W}_{\ell}^{\ast} \otimes {\bigwedge}  ^{i} (\mathfrak
{p}^{-}) ^{\ast} \otimes
\calS(V^n) \otimes {\bigwedge} ^{j}\mathfrak{p^{\ast}} \otimes  T^{\ell}(V) 
\right]^{K' \times K},
\]
where the multiplication $\cdot$ in $C$ is given componentwise. For each
$\ell$, we have a double complex $(C^{\bullet,\bullet}_{\ell},
\bar{\partial},d)$ with commuting differentials
\begin{align*}
\bar{\partial} & = \sum_{j=1}^{n(n+1)/2} 1 \otimes A(\bar{\eta}_j) \otimes 
\omega(\bar{Z_j}) \otimes 1\otimes 1, \\
d  & = d_{\calS} + d_{V},
\end{align*}
where
\begin{align*}
d_{\calS} & = \sum_{\substack{\alpha,\mu}} 1 \otimes 1 \otimes \omega(X_
{\alpha\mu}) \otimes A(\omega_{\alpha\mu}) \otimes  1, \\
d_V & = \sum_{\substack{\alpha,\mu}} 1 \otimes 1 \otimes 1 \otimes A(\omega_
{\alpha\mu}) \otimes \rho(X_{\alpha\mu}).
\end{align*}
Here $A(\cdot)$ denotes the left multiplication, while $\rho$ is the
derivation action of $\mathfrak{g}$ on $T^{\ell}(V)$. Furthermore, $K'$ acts
on the first three tensor factors of $C^{i,j}_{\ell}$, while $K$ acts
on the last one. The actions on $S(V^n)$ are given by the Weil representation,
while the actions on the other tensor factors are the natural ones.

We also have an analogous complex $C^{\bullet,\bullet}_{[\la]}$ by replacing $W^
{\ast}_{\ell}$ and
$T^{\ell}(V)$ with $W^{\ast}_{\la}$ and $S_{[\la]}(V)$ respectively.

We call a $d$-closed element $\varphi \in C^{i,j}$  {\em holomorphic} if the
cohomology class $[\varphi]$ is $\bar{\partial}$-closed, i.e., there
exists $\psi \in  C^{i+1,j-1}$ such that
\[
\bar{\partial} \varphi = d \psi.
\]

Note that the maps $\bar{\partial}$ and $d$ correspond to the usual operators 
$\bar{\partial}$ and $d$ under the isomorphism
\[
[ {W}^{\ast}_{\ell} \otimes \mathcal{A}^{0,i}(\h_n) \otimes \mathcal{S}(V^n) 
\otimes \mathcal{A}^j(D) \otimes  T^{\ell}(V)]^{G' \times G} \to C^{i,j}_{\ell},
\]
given by evaluation at the base points of $\h_n$ and $D$. In the
following we will frequently identify these two spaces, and by
abuse of notation we will use the same symbol for corresponding
objects.

\subsection{Special Schwartz forms}

We construct for $n \leq p$ a family of Schwartz functions
$\varphi_{nq,\ell}$ on $V^n$ taking values in $\mathcal{A}^{nq}(D)
\otimes W^{\ast}_{\ell} \otimes T^{\ell}(V)$, the space of differential
$nq$-forms on $D$ which take values in $   {W^{\ast}_{\ell}
  \otimes  T^{\ell}(V)}$. That is, $\varphi_{nq,\ell} \in  C^{0,j}_{\ell}$:
\begin{align*}\label{varphi}
\varphi_{nq,\ell} &\in \left[W_{\ell}^{\ast} \otimes \calS(V^n) \otimes \mathcal
{A}^{nq}(D) \otimes  T^{\ell}(V)\right]^{K'\times G} \\
& \notag \simeq \left[  W_{\ell}^{\ast}  \otimes
\calS(V^n) \otimes {\bigwedge} ^{nq}(\mathfrak{p^{\ast}}) \otimes
T^{\ell}(V) \right]^{K' \times K}.
\end{align*}
These Schwartz forms are the generalization of the 'scalar-valued' Schwartz 
forms considered by Kudla and Millson \cite{KMI,KMII,KM90} to the coefficient 
case. 



Starting with the standard Gaussian,
\[
\varphi_0(\x) = e^{-\pi tr(x,x)_{z_0}} \in \mathcal{S}(V^n),
\]
with $\x =(x_1, \cdots, x_n) \in V^n$, the `scalar-valued' form $\varphi_{nq,0}
$ is given by applying the operator
\begin{gather*}
\calD:  \mathcal{S}(V^n) \otimes {\bigwedge} ^{\ast} (\mathfrak{p^{\ast}})
\longrightarrow   \mathcal{S}(V^n) \otimes {\bigwedge} ^{\ast + nq}(\mathfrak{p^
{\ast}}), \\
\calD = \frac{1}{2^{nq/2}} \prod_{i=1}^n \prod_{\mu = p+1}^{p+q} \left[
\sum_{\alpha =1}^{p} \left(  x_{\alpha i} - \frac{1}{2\pi}
\frac{\partial}{\partial x_{\alpha i }}  \right)   \otimes
A(\omega_{\alpha\mu}) \right],
\end{gather*}
to  $\varphi_0 \otimes 1  \in   [S(V^n) \otimes {\bigwedge} ^0(\mathfrak{p^
{\ast}}) ]^K$:
\begin{equation*}
 \varphi_{nq,0} = \calD (\varphi_0 \otimes 1 ).
\end{equation*}
Note that this is $2^{nq/2}$ times the corresponding quantity in \cite{KM90}. 
We have
\[
\varphi_{nq,0}(\x) \in C_0^{0,nq}= \left[  \C(\chi_{-m/2}) \otimes    \calS
(V^n) \otimes
{\bigwedge} ^{nq}(\mathfrak{p^{\ast}}) \right ]^{K' \times K}
\]
Here the $K$-invariance is immediate, while the $K'$-invariance is Theorem~3.1 
in \cite{KMI}.

We let
\[
\mathcal{A} = End_{\C} \left( \calS(V^n) \otimes {\bigwedge} ^{\ast}(\mathfrak
{p^{\ast}}) \otimes T(V) \right),
\]
where $T(V) = \bigoplus_{\ell=0}^{\infty} T^{\ell}(V)$ denotes the 
(complexified) tensor algebra of $V$. Note that $\mathcal{A}$ is an associative 
$\C$-algebra by composition. We now define for $1\leq i \leq n$ another 
differential operator $\calD_i \in \mathcal{A}$ by
\[
\calD_i = \frac12 \sum_{\alpha=1}^p \left(   x_{\alpha i} - \frac{1}{2\pi}
\frac{\partial}{\partial x_{\alpha i}}  \right)   \otimes 1 \otimes 
A(e_{\alpha}).
\]
Here $A(e_{\alpha})$ denotes the left multiplication by $e_{\alpha}$ in $T(V)$. 
Note that the operator $\calD_i$ is clearly $K$-invariant.
We introduce a homomorphism $T: \C^n \mapsto \mathcal{A}$ by
\[
T(\eps_i) = \calD_i,
\]
where $\eps_1,\dots,\eps_n$ denotes the standard basis of $\C^n$. Let $m_
{\ell}: T^{\ell} \calA \mapsto \calA$ be the $\ell$-fold multiplication. We now 
define 
\[
\calT_{\ell}: T^{\ell} (\C^n) \longrightarrow \calA
\]
by
\[
\calT_{\ell} = m_{\ell} \circ \left( \bigotimes{} ^{\ell} T \right).
\]

We identify 
\[
\Hom_{\C}\left( W_{\ell}, \mathcal{S}(V^n) \otimes {\bigwedge}
  ^{nq}(\mathfrak{p^{\ast}}) \otimes T^{\ell}(V) \right) \simeq 
 W_{\ell}^{\ast} \otimes \mathcal{S}(V^n) \otimes {\bigwedge}
  ^{nq}(\mathfrak{p^{\ast}}) \otimes T^{\ell}(V),
\]
and use the same symbols for corresponding objects. 

\begin{definition}
We define 
\[
\varphi_{qn,\ell} \in \Hom_{\C}\left( W_{\ell}, \mathcal{S}(V^n) \otimes 
{\bigwedge} ^{nq}(\mathfrak{p^{\ast}}) \otimes T^{\ell}(V) \right)^K
\]
by
\[
\varphi_{nq,\ell}(w) = \calT_{\ell}(w) \varphi_{nq,0}
\]
for $w \in T^{\ell} (\C^n)$. We put $ \varphi_{nq,\ell} =0$ for $\ell <0$. 
\end{definition}


\smallskip

Note that the symmetric group $S_{\ell}$ on $\ell$ letters is acting
on $T^{\ell}(\C^n)$ and $T^{\ell} (V)$ in the natural fashion, which
gives rise to a natural action of $S_{\ell}$ on $C^{i,j}_{\ell}$. We will
now show that $\varphi_{nq,\ell}$ is an equivariant map with respect
to ${S}_{\ell}$. More precisely, we have

\begin{proposition}

\[
\varphi_{nq,\ell}  \in \Hom_{\C}\left( W_{\ell}, \mathcal{S}(V^n)
  \otimes {\bigwedge} ^{nq}(\mathfrak{p^{\ast}}) \otimes T^{\ell}(V)
\right)^{S_{\ell}\times K},
\]
that is,
\[
\varphi_{nq,\ell} \circ s = (1 \otimes 1 \otimes s) \varphi_{nq,\ell}.
\]

\end{proposition}

\begin{proof}

We first need

\begin{lemma}
Let $s \in S_{\ell}$. Then
\[
\calT_{\ell} \circ s = (1 \otimes 1 \otimes s) \circ \calT_{\ell}.
\]
\end{lemma}

\begin{proof}
Let $i_1,i_2,\dots,i_{\ell} \in \{1,\dots,n\}$. Then
\[
\left( \bigotimes{} ^{\ell} T \right) \left( s (\eps_{i_1} \otimes \cdots 
\otimes \eps_{i_{\ell}})\right) = s \left(T(\eps_{i_1}) \otimes \cdots \otimes T
(\eps_{i_{\ell}})   \right),
\]
where $s$ on the right hand side permutes the factors of $\calA \otimes \cdots 
\otimes \calA$. Hence
\begin{align*}
\calT_{\ell} \left(   s (\eps_{i_1} \otimes \cdots \otimes \eps_{i_{\ell}})    
\right) = m_{\ell} \left(s\left(T(\eps_{i_1}) \otimes \cdots \otimes T(\eps_{i_
{\ell}})   \right)\right), 
\end{align*}
thus $\calT_{\ell} \left(   s (\eps_{i_1} \otimes \cdots \otimes \eps_{i_
{\ell}})    \right)$ takes the factors of 
$\left( \bigotimes{} ^{\ell} T \right) \left( s (\eps_{i_1} \otimes \cdots 
\otimes \eps_{i_{\ell}})\right)$, permutes them according to $s$ and then 
multiplies them in $\calA$. The product $\calT_{\ell} (\eps_{i_1} \otimes 
\cdots \otimes \eps_{i_{\ell}})$ is a sum of tensor products of products of 
certain differential operators in $\calS(V^n)$ with products of the $e_{\alpha}
$'s in $T(V)$. But the differential operators in  $\calS(V^n)$ commute with 
each other so that the rearrangement of the $D_i$'s has no effect on this 
factor in the tensor product. Hence $s$ only acts on the third factor of  
$\mathcal{S}(V^n) \otimes {\bigwedge} ^{nq}(\mathfrak{p^{\ast}}) \otimes T^
{\ell}(V)$.
\end{proof}

The proposition now follows easily. For $w \in T^{\ell} (\C^n$), we have
\begin{align*}
\varphi_{nq,\ell} ( s w) = \calT_{\ell}(s w) \varphi_{nq,0} = \left(  (1 
\otimes 1 \otimes s)  \calT_{\ell}(w) \right) \varphi_{nq,0} =   (1 \otimes 1 
\otimes s) \left(    \calT_{\ell}(w) \varphi_{nq,0} \right).
\end{align*}
\end{proof}

This will now enable us to define the Schwartz forms
$\varphi_{nq,[\la]}$. We first note

\begin{lemma}
For any standard filling $t(\la)$ of $D(\la)$, the composition
\[
(1 \otimes 1 \otimes \pi_{t(\la)}) \circ \varphi_{nq,\ell} : T^{\ell} \C^n 
\longrightarrow 
 \mathcal{S}(V^n) \otimes {\bigwedge} ^{nq}(\mathfrak{p^{\ast}}) \otimes S_{t
(\la)}(V)
\]
descends to a map
\[
S_{t(\la)} (\C^n) \longrightarrow  \mathcal{S}(V^n) \otimes {\bigwedge} ^{nq}
(\mathfrak{p^{\ast}}) \otimes S_{t(\la)}(V).
\]
\end{lemma}

\begin{proof}
We have $(1 \otimes 1 \otimes s_{f(\la)})^2 = 1 \otimes 1 \otimes s_{f(\la)}$. 
Since $\varphi_{nq,\ell}$ is equivariant with respect to $S_{\ell}$, we have 
\[
(1 \otimes 1 \otimes s_{f(\la)}) \circ \varphi_{nq,\ell} ( w_1 \otimes \cdots 
\otimes w_{\ell}) = 
(1 \otimes 1 \otimes s_{f(\la)}) \circ \varphi_{nq,\ell} (s_{f(\la)} (w_1 
\otimes \cdots \otimes w_{\ell}))
\]
for all $w_i \in \C^n, 1 \leq i \leq \ell$.
\end{proof}

We use the lemma for the standard filling $t_0(\la)$ to introduce $\varphi_{nq,
[\la]}$. 

\begin{definition}
We define
\[
\varphi_{nq,[\la]} \in \Hom_{\C}\left( S_{\la}(\C^n), \mathcal{S}(V^n) \otimes 
{\bigwedge} ^{nq}(\mathfrak{p^{\ast}}) \otimes S_{[\la]}(V) \right)^K
\]
by
\[
\varphi_{nq,[\la]}(w) = (1 \otimes 1 \otimes \pi_{[\la]})(\varphi_{nq,\ell}
(\iota_{\la}(w)),
\]
the projection onto  $S_{[\la]}(V)$, the harmonic tensors in $S_{\la}(V)$. 

\end{definition}

\subsection{Fundamental Properties of the Schwartz forms}

We will now state the four basic properties of our Schwartz
forms. These are:
\begin{itemize}
\item
$K'$-invariance; thus $\varphi_{nq,\ell} \in C^{0,nq}_{\ell}$
\item
$d$-closedness; thus  $\varphi_{nq,\ell}$ defines a cohomology class  $[\varphi_
{nq,\ell}]$
\item
The holomorphicity of  $[\varphi_{nq,[\la]}]$
\item
A recursion formula relating $[\varphi_{nq,\ell}]$ to $[\varphi_{nq,\ell-1}]$
\end{itemize}
The first three properties are the generalizations of the properties
of $\varphi_{nq,0}$ in \cite{KMI,KMII,KM90}, the trivial coefficient
case. Except for the $K'$-invariance, we will reduce the statements to
the case of $n=1$. Our main tool in proving these properties will be
then the Fock model of the Weil representation. We will carry out the
proofs for the $K'$-invariance and for the other statements in the
case of $n=1$ in the next section.

\begin{theorem}\label{MAIN3}

The forms $\varphi_{nq,\ell}$ and 
$\varphi_{nq,[\la]}$ are $K'$-invariant, i.e., 
\[
\varphi_{nq,\ell} \in C_{\ell}^{0,nq} = 
  [  W_{\ell}^{\ast} \otimes
\calS(V^n) \otimes {\bigwedge} ^{nq}(\mathfrak{p^{\ast}}) \otimes  T^{\ell}(V) ]
^{K'\times K}
\]
and
\[
\varphi_{nq,[\la]}  \in C_{[\la]}^{0,nq} = 
  [ W_{\lambda^{\prime}}^{\ast} \otimes
\calS(V^n) \otimes {\bigwedge} ^{nq}(\mathfrak{p^{\ast}}) \otimes  S_{[\la]}
(V) ]^{K'\times K}.
\]
In particular, for $n=1$, we have
\[
\varphi_{q,\ell} \in \left[ \C(\chi_{-\ell-m/2}) \otimes \calS(V) \otimes 
{\bigwedge} ^{q}(\mathfrak{p^{\ast}}) \otimes  T^{\ell}(V) \right]^{K \times 
K'}.
\]

\end{theorem}

\begin{proof}
We consider the first statement using the Fock model in the next section. The 
second statement follows from the first by projecting onto $S_{[\la]}(V)$.
\end{proof}

The $K'$-invariance of the Schwartz forms will enable us in
Section~\ref{Theta} to construct theta series using the forms $\varphi_{nq,
[\la]}$.

\begin{theorem}\label{MAIN1}

The forms $\varphi_{nq,\ell}$ and $\varphi_{nq,[\la]}$ define closed 
differential forms on $D$, i.e.,
\[
d \varphi_{nq,\ell}(\x) =0
\]
for all $\x \in V^n$. In particular, $\varphi_{nq,[\la]}(\x)$ defines a 
(deRham) cohomology class
\[
[\varphi_{nq,[\la]}(\x)] \in H^{nq}\left(D, {\Hom}_{\C}(S_{\la}(\C^n), S_{[\la]}
(V)) \right).
\]

\end{theorem}

\begin{proof}
We will prove the case $n=1$ in the next section using the Fock model of the 
Weil representation. For general $n$, it is enough to show that $\varphi_
{nq,\ell}(\eps_{i_1} \otimes \cdots \otimes \eps_{i_{\ell}})$ is closed for any 
$n$-tuple $(\eps_{i_1}, \dots, \eps_{i_{\ell}})$. By the $S_{\ell}$-
equivariance of $\varphi_{nq,\ell}$ we can assume that ${i_1} \leq \cdots \leq 
\ i_{\ell}$, so that 
$\eps_{i_1} \otimes \cdots \otimes \eps_{i_{\ell}} = \eps_1^{\otimes \ell_1} 
\otimes \cdots \otimes \eps_n^{\otimes \ell_n}$ for some non-negative integers 
$\ell_1,\dots,\ell_n$. But this implies that
\[
\varphi_{nq,\ell}(\eps_{1}^{\otimes \ell_1} \otimes \cdots \otimes \eps_{n}^
{\otimes \ell_n})(\x) = \varphi_{q,\ell_1}(x_1) \wedge \cdots \wedge \varphi_
{q,\ell_n}(x_n).
\]
Here the wedge $\wedge$ means the usual wedge for $\mathcal{A}(D)$ and the 
tensor product in the other slots. This reduces the closedness of $\varphi_
{nq,\ell}$ to the case $n=1$.
\end{proof}

To state the last two properties of the forms $\varphi_{nq,[\la]}$, we
first need to introduce some more notation. We define a map
\[
\sigma: \C^n \longrightarrow \left(V^{n}\right)^{\ast} \otimes {\bigwedge} ^
{\ast} \mathfrak{p}^{\ast} \otimes V
\]
by 
\[
\sigma(\epsilon_i) = \sum_{j=1}^m x_{ij}  \otimes 1 \otimes e_j.
\]
Here the $x_{ij}, 1 \leq i \leq n, 1 \leq j \leq m$ are the standard 
coordinates on $V^n$.
Thus under the identification of $\left(V^{n}\right)^{\ast} \otimes {\bigwedge} 
^{\ast} \mathfrak{p}^{\ast} \otimes V$
with $Hom(V^n, {\bigwedge} ^{\ast} \mathfrak{p}^{\ast} \otimes V)$ we have 
\begin{align*}
\sigma (\eps_i)(\x) = \sum_{j=1}^m x_{ij} ( 1 \otimes e_j) 
& =  1 \otimes \sum_{j=1}^m x_{ij} e_j 
 =  1\otimes x_i. 
\end{align*}
Now the $x_{ij}$ are numbers, the coordinates of the $n$-tuple of vectors $\x$.
Under the identification of $\mathfrak{p}^{\ast} \otimes V$ with $1 \otimes 
\mathfrak{p}^{\ast} \otimes V$ we may
rewrite the above formula as
$$\sigma (\eps_i)(\x) = 1 \otimes 1\otimes x_i.$$
By interpreting $\sigma(\eps_i)$ as the left multiplication operator by $\sigma
(\eps_i)$ we can interpret $\sigma$ as a map from $\C^n$ to 
$\mathcal{A}$ (and we do not distinguish between these two interpretations). We 
let $\sigma_{\ell}$ be the $\ell$-th (exterior) tensor 
power of $\sigma$, and for $\lambda$ a partition of $(\ell)$, we put
\[
\sigma_{\lambda} = \sigma_{\ell} \circ \iota_{\lambda}:  S_{\lambda} \C^n 
\longrightarrow \mathcal{A}.
\]
Note that we do not need to distinguish between $\lambda^{\prime}$ and 
$\lambda$ because only the nonzero parts of the partition
matter here.

\begin{lemma}\label{auxlemma}

\begin{itemize}

\item[(i)]
Let $f$ be a semistandard filling of $D(\la$). Then
\[
\sigma_{\la} (\eps_{f(\la)})(\x) = 1 \otimes 1 \otimes \x_{f(\la)}
\]
for any $\eps = \eps_{i_1} \otimes \cdots \otimes \eps_{i_{\ell}} \in
T^{\ell}(\C^n)$. In particular,
\[
\sigma_{\la}: S_{\la} \C^n \longrightarrow \left(V^{n}\right)^{\ast}
\otimes {\bigwedge} ^{\ast} \mathfrak{p}^{\ast} \otimes S_{\la}(V).
\]

\item[(ii)]
The map $\sigma_{\la}$ is $GL_n(\C)$-invariant, i.e.,
\[
\sigma((a^{-1} \eps)_{f(\la)}) (\x a) = \sigma(\eps_{f(\la)}) (\x)
\]
for $a \in GL_n(\C)$.
\end{itemize}

\end{lemma}

\begin{proof}
For (i), first note
\[
\sigma_{\ell} (\eps_{i_1} \otimes \cdots \otimes \eps_{i_{\ell}})(\x) =
1 \otimes 1 \otimes (x_{i_1} \otimes \cdots \otimes x_{i_{\ell}}).
\]
Indeed,
\begin{align*}
\sigma_{\ell} (\eps_{i_1} \otimes \cdots \otimes \eps_{i_{\ell}})(\x) = \sigma_
{\ell} (\eps_{i_1}) \circ \cdots \circ 
\sigma_{\ell} (\eps_{i_{\ell}}) (\x)  &=
(1 \otimes 1 \otimes x_{i_1}) \circ \cdots \circ (1 \otimes 1 \otimes x_{i_
{\ell}}) \\
&= 1 \otimes 1 \otimes (x_{i_1} \otimes \cdots \otimes x_{i_{\ell}}).
\end{align*}
But now for $s \in S_{\ell}$ and $w_i \in \C^n$, $1 \leq i \leq \ell$, we have
\[
\sigma_{\ell} \left(s (w_1 \otimes \cdots \otimes w_{\ell})  \right) = (1 
\otimes 1 \otimes s) \sigma(w_1 \otimes \cdots \otimes w_{\ell}),
\]
which gives immediately
\[
\sigma_{\ell}(\eps_{f(\la)})(\x) = 1 \otimes 1 \otimes \x_{f(\la)},
\]
as claimed. (ii) follows easily from $ \sigma( a^{-1} \eps)(\x a) = \sigma(\eps)
(\x)$.
\end{proof}

We can therefore define $\sigma_{[\la]}$ by postcomposing with the
harmonic projection $\pi_{[\la]}$ onto $S_{[\la]}(V)$, and we have
\[
\sigma_{[\la]}(\pi_{\la}\eps_{f(\la)})(\x) = 1 \otimes 1 \otimes \x_{[f(\la)]}.
\]

\smallskip

For $v \in V$, we let $A_j(v):T^{\ell-1}(V) \to T^{\ell}(V)$ be the insertion 
of $v$ into the $j$-th spot. We let $A_{jk}: T^{\ell-2}(V) \to T^{\ell}$
\[
A_{jk}(f) = \sum_{\alpha=1}^p A_j(e_{\alpha}) A_k(e_{\alpha}) - \sum_{\mu =p+1}^
{p+q} A_j(e_{\mu}) A_k(e_{\mu})
\]
be the insertion of the invariant metric into the $(j,k)$-th spot, and we put
\[
A(f) =  \frac{1}{2} \sum_{j=1}^{\ell} \sum_{k=1}^{\ell-1}A_{jk}(f).
\]

One of the fundamental properties of the scalar-valued Schwartz form
$\varphi_{nq,0}$ is that for $(\x,\x)$ positive semidefinite,
$\varphi_{nq,0}(\x)$ gives rise to a Thom form for the special cycle $C_{X}$. 
In view of Lemma~{\ref{PDform}}, we now relate $\varphi_{nq,[\la]}(\x)$ to
$\sigma_{[\la]}(\x)\varphi_{nq,0}(\x)$.

\begin{theorem}\label{MAIN2}

\begin{itemize}

\item[(i)]
Let $n=1$ and let $\sigma_j$ be the operator on  $\mathcal{S}(V) \otimes 
{\bigwedge} ^{\ast} \mathfrak{p}^{\ast} \otimes T(V)$ defined by $\sigma_j(x) = 
1 \otimes 1 \otimes A_j(x)$. Then for each $j =1,\dots,\ell$, we have in 
cohomology
\[
 [\varphi_{q,\ell}]  = [\sigma_j \varphi_{q,\ell-1}] \, + \, \frac{1}{4\pi} 
\sum_{k=1}^{\ell-1} [A_{jk}(f)\varphi_{q,\ell-2}]
\]
for all $x \in V$. In particular,
\
\[
 [\varphi_{q,[\ell]}]  = [\sigma_{[\ell]} \varphi_{q,0}].
\]

\item[(ii)]
For general $n$, we have in cohomology
\[
[\varphi_{nq,[\la]}] = [\sigma_{[\la]}\varphi_{nq,0}],
\] 
i.e,, for all semistandard fillings $f$,
\[
[\varphi_{nq,[\la]}(\eps_{f(\la)}) (\x)] = \left[(1 \otimes 1 \otimes
\x_{[f(\la)]}) \varphi_{nq,0}(\x)\right],
\]
where $\eps = \eps_{i_1} \otimes \cdots \otimes \eps_{i_{\ell}} \in T^{\ell}
(\C^n)$.

\end{itemize}
\end{theorem}

\begin{proof}
We proof (i) in the next section via the Fock model. For (ii), we first see by 
(i) that up to exact forms we have
\begin{align*}
\varphi_{nq,\ell} (\eps_1^{\otimes \ell_1} \otimes \cdots \otimes \eps_n^
{\otimes \ell_n})(\x) & = 
\left( \sigma(x_1) \varphi_{q,\ell_1-1}(x_1) +   \frac{1}{4\pi} \sum_{k=1}^
{\ell_1-1} [A_{jk}(f) \varphi_{q,\ell_1-1}(x_1) \right) \\
& \quad  \wedge \cdots \\
& \quad  \wedge  \left( \sigma(x_n) \varphi_{q,\ell_{n-1}}(x_n) 
+   \frac{1}{4\pi} \sum_{k=1}^{\ell_n-1} [A_{jk}(f) \varphi_{q,\ell_n-1}(x_n)
\right),
\end{align*}
and we get a similar statement for $\varphi_{nq,\ell}(\eps_{i_1}
\otimes \cdots \otimes \eps_{i_{\ell}})$ by the $S_{\ell}$-equivariance of
$\varphi_{nq,\ell}$. Iterating and using Lemma~\ref{auxlemma} then 
gives a statement for $\varphi_{nq,\ell}(\x)$ 
with coefficients in $S_{\la}(V)$ analogous to $(ii)$ - up to terms coming from 
the metric. Projecting to $S_{[\la]}(V)$ now gives the claim.
\end{proof}

One of the main results of \cite{KM90} is that the scalar-valued
cohomology class $[\varphi_{nq,0}]$ is holomorphic, i.e, 
$[\overline{\partial}\varphi_{nq,0}]=0$. We will now show that the
more general cohomology classes $[\varphi_{nq,[\la]}(\x)]$ are
holomorphic as well.

For $n=1$, we have $\mathfrak{g}' = \mathfrak{sl}_2(\C)$, and the anti-
holomorphic tangent space $\mathfrak{p}^-$ is spanned by the element $L =  
\tfrac12 \left( \begin{smallmatrix}1 & -i \\ -i & -1 \end{smallmatrix} \right)
$. The Weil representation action of $L$ corresponds to the classical Maass 
lowering operator  $-2i v^2 \tfrac{\partial}{\partial \bar{\tau}}$ on the upper 
half plane.

\begin{theorem}\label{MAIN4}

\begin{itemize}

\item[(i)]
Let $n=1$.
Then in cohomology, we have
\[
[\omega(L)\varphi_{q,\ell}] =\frac{-1}{4 \pi} [A(f)\varphi_{q,\ell-2}]. 
\]
In particular,
\[
[\overline{\partial} \varphi_{q,[\ell]}] = 0.
\]

\item[(ii)]
For general $n$, we have 
\[
[\overline{\partial} \varphi_{nq,[\la]}] =0.
\]


\end{itemize}

\end{theorem}

\begin{proof}
We will prove (i) in the next section. (ii) follows from (i) by generalizing 
the argument given for the scalar valued case in \cite{KM90}, Theorem~5.2. 
First note that we have to show $\overline[{\partial}_{ij} \varphi_{nq,[\la]}]=0
$ for all $n(n+1)/2$ partial derivatives $\overline{\partial}_{ij}$ ( $i \leq 
j$) in $\mathfrak{p}^-$. By (i), we see
\[
(1 \otimes 1 \otimes \pi_{[\la]})
\overline{\partial_{ii}} \varphi_{nq,\ell}(\eps_1^{\otimes \ell_1} \otimes 
\cdots \otimes \eps_{n}^{\otimes \ell_n}) = 0
\]
up to an exact form. By the $S_{\ell}$-equivariance, we then see
\[
(1 \otimes 1 \otimes \pi_{[\la]})
\overline{\partial_{ii}} \varphi_{nq,\ell}(\eps_{i_1} \otimes \cdots \otimes 
\eps_{i_{\ell}}) = 0,
\]
again, up to an exact form. This gives the desired vanishing for the anti-
holomorphic tangent space $\mathfrak{p}_0^-$ of $\h \times \cdots \times \h$, 
naturally embedded into $\h_n$. By the $K'$-invariance of $\varphi_{nq,\ell}$ 
we now see that $(1 \otimes 1 \otimes c_{[\la]}) \varphi_{nq,\ell}$ is 
annihilated by the $Ad \, K'$ orbit of $\mathfrak{p}_0^-$ inside $\mathfrak{p}^-
$, which is all of  $\mathfrak{p}^-$.
\end{proof}

\section{Proof of the fundamental properties of the Schwartz forms}\label{Fock}

The purpose of this section is to prove the $K'$-invariance of 
$\varphi_{nq,\ell}$ and for $n=1$ the other fundamental properties of $\varphi_
{q,\ell}$.
given in the previous section. Our main tool will be the Fock model of the Weil 
representation, which we review in the appendix.

By abuse of notation we will frequently use in the following the same symbols 
for corresponding objects and operators in the two models.

\subsection{The Schwartz forms in the Fock model and the $K'$-invariance}


For multi-indices $\underline{\alpha} = (\alpha_1,\cdots,\alpha_q)$ and 
$\underline{\beta} = (\beta_1,\cdots,\beta_{\ell})$, (usually suppressing their 
length), we will write
\begin{gather*}
\omega_{\underline{\alpha}} = \omega_{\alpha_1 p+1} \wedge \cdots \wedge \omega_
{\alpha_q p+q}, \\
z_{\underline{\alpha}j} = z_{\alpha_1j} \cdots z_{\alpha_qj}, \\
e_{\underline{\beta}} = e_{\beta_1} \otimes \cdots \otimes e_{\beta_{\ell}}.
\end{gather*}
 Here we have returned to our original notation, denoting the standard basis 
elements of $V$ by $e_{\alpha}$ and $e_{\mu}$. In the Fock model, the ``scalar-
valued'' Schwartz form $\varphi_{nq,0}$ becomes with this notation
\[
\varphi_{nq,0} = \frac{1}{2^{nq/2}} \left(\frac{-i}{2\pi}\right)^{nq}  \sum_
{{\underline{\alpha_1}},\dots, {\underline{\alpha_n}}} z_{{\underline{\alpha_1}
1}} \cdots z_{{\underline{\alpha_n}n}} \otimes \omega_{{\underline{\alpha_1}1}} 
\wedge \cdots \wedge \omega_{{\underline{\alpha_n}n}} \otimes 1 \in \mathcal{F} 
\otimes {\bigwedge} ^{nq}(\mathfrak{p}^{\ast}) \otimes T^0(V)
\]
We define
\[
\varphi_{0,\ell} \in \Hom_{\C}(T^{\ell}(\C^n), \mathcal{F}\otimes {\bigwedge} ^
{0}(\mathfrak{p}^{\ast})  \otimes  T^{\ell}(V))
\]
by 
\[
\varphi_{0,\ell}(\eps_{i_1} \otimes \cdots \otimes \eps_{i_{\ell}}) = 
\left( \tfrac{-i}{4\pi}\right)^{n\ell} \sum_{\underline{\beta}} z_{\beta_1 
i_1 } \cdots 
z_{\beta_{\ell i_{\ell}}} \otimes  1 \otimes e_{\underline{\beta}} .
\]
We then easily see

\begin{lemma}
\[
\varphi_{nq,\ell} = \varphi_{nq,0} \cdot \varphi_{0,\ell},
\]
where the multiplication is the natural one in $ \Hom_{\C}(T(\C^n),\mathcal{F} 
\otimes {\bigwedge} ^{\ast}(\mathfrak{p}^{\ast})  \otimes T(V))$.

\end{lemma}

We should note that only in the Fock model we have such a ``splitting'' of 
$\varphi_{nq,\ell}$ into the product of two elements. We do not have an 
analogous statement in the Schroedinger model (only in terms of operators 
acting on the Gaussian $\varphi_0$).


\begin{theorem}[Theorem~\ref{MAIN3}]\label{FockK'}
The form $\varphi_{nq,\ell}$ is $K'$-invariant.
\end{theorem}

\begin{proof}
We show this on the Lie algebra level. The element $k' = \frac{1}{2i} w'_j 
\circ w_k''  \in \mathfrak{k} \simeq \mathfrak{gl}_n(\C)$ is the endomorphism 
of $\C^n$ mapping $\eps_j$ to $\eps_k$ and annihilating the other basis 
elements. To show $\omega(k') \varphi_{nq,\ell} =0$, we need to show
\[
\omega(k') \left( \varphi_{nq,\ell}(\eps_{i_1} \otimes \cdots \otimes \eps_{i_
{\ell}}) \right) = 
\varphi_{nq,\ell}( k'( \eps_{1_1} \otimes \cdots \otimes \eps_{1_n})).
\]

>From Lemma~\ref{Fock2} we see
\begin{align*}
\omega(k') \left( \varphi_{nq,\ell}(\eps_{i_1} \otimes \cdots \otimes \eps_{i_
{\ell}}) \right) &= 
\omega(k') \varphi_{nq,0} \cdot \varphi_{0,\ell}(\eps_{1_1} \otimes \cdots 
\otimes \eps_{1_n}) \\
& \quad + \varphi_{nq,0} \cdot \sum_{\alpha =1}^p z_{\alpha k} \frac{\partial}
{\partial z_{\alpha j}}
\left( \varphi_{0,\ell}(\eps_{1_1} \otimes \cdots \otimes \eps_{1_n}) \right).
 \end{align*}
We have $ \omega(k') \varphi_{nq,0} = 0 $, since $
 \varphi_{nq,0} \in  \left[\C(\chi_{-m/2}) \otimes \mathcal{F} \otimes  
\bigwedge^{nq} (\mathfrak{p}^{\ast} )\right]^{K'}$ by \cite{KMI}, Theorem~5.1. 
On the other hand, one easily sees
\[
 \sum_{\alpha =1}^p z_{\alpha k} \frac{\partial}{\partial z_{\alpha j}}
\left( \varphi_{0,\ell}(\eps_{1_1} \otimes \cdots \otimes \eps_{1_n}) \right) = 
\varphi_{0,\ell}( k'( \eps_{1_1} \otimes \cdots \otimes \eps_{1_n})).
\]
The assertion follows. 
\end{proof}

\subsection{The Schwartz forms for $n=1$}

For $n=1$, we consider the forms $\varphi_{q,\ell}$,  $\varphi_{q,0}$, and 
$\varphi_{0,\ell}$ to be in  $\calF \otimes {\bigwedge} {^{\ast}}(\mathfrak{p}^
{\ast})\otimes T(V)$, and we have
\begin{align*}\label{Phi}
\varphi_{q,\ell} 
=  c_{q,\ell} \sum_{\underline{\alpha},\underline{\beta}} z_{\underline
{\alpha}} z_{\underline{\beta}} \otimes \omega_{\underline{\alpha}} \otimes e_
{\underline{\beta}}.
\end{align*}
Here $c_{q,\ell} ={2^{q/2}}(-i/4\pi)^{q+\ell}$. Also
\begin{gather*}
\varphi_{q,0}= 2^{q/2}(-i/4\pi)^{q} 
\sum_{\underline{\alpha}} z_{\underline{\alpha}} \otimes 
\omega_{\underline{\alpha}} \otimes 1 \qquad \text{and} \qquad 
\varphi_{0,\ell} = (-i/4\pi)^{\ell} \sum_{\underline{\beta}} z_{\underline
{\beta}} \otimes 1 \otimes  e_{\underline{\beta}}.
\end{gather*}

For later use, we note that Theorem~\ref{FockK'} for $n=1$ boils down to
\begin{equation}\label{K'1}
 \sum_{\alpha =1}^p \left(z_{\alpha} \frac{\partial}{\partial z_{\alpha }}
\otimes 1 \otimes 1 \right) \varphi_{q,\ell} = (q+\ell) \, \varphi_{q,\ell},
\end{equation}
which follows directly from 
\begin{equation}\label{K'2}
 \sum_{\alpha =1}^p z_{\alpha} \frac{\partial}{\partial z_{\alpha }}
\varphi_{q,0} = q \, \varphi_{q,0} \qquad \text{and} \qquad 
 \sum_{\alpha =1}^p z_{\alpha} \frac{\partial}{\partial z_{\alpha }}
\varphi_{0,\ell} = \ell \, \varphi_{0,\ell}.
\end{equation}


\subsection{Closedness}

Similarly to  the Schroedinger model, the differentiation $d$ in the Lie 
algebra complex $
\mathcal{F} \otimes \bigwedge{^{\ast}}(\mathfrak{p}^{\ast}) \otimes S^{\ell}(V)$
is given by $d = d_{\calF} + d_V$ with 
\begin{gather}\label{d1}
d_{\calF} = \sum_{\alpha,\mu} \omega (X_{\alpha\mu}) \otimes
A(\omega_{\alpha\mu})\otimes 1 \qquad \text{and} \qquad
d_V = \sum_{\alpha,\mu} 1 \otimes 
A(\omega_{\alpha\mu})\otimes \rho(X_{\alpha\mu}).
\end{gather}
Furthermore, we write $d_{\calF} = d_{\calF}'+ d_{\calF}''$ with 
\begin{gather}\label{dF}
d_{\calF}' = -4 \pi \sum_{\alpha,\mu}   \frac{\partial^2}{\partial
  z_{\alpha } \partial z_{\mu }}  \otimes
A(\omega_{\alpha\mu})\otimes 1, \\
d_{\calF}''=   \frac{1}{4\pi} \sum_{\alpha,\mu} z_{\alpha }z_{\mu } \otimes
A(\omega_{\alpha\mu})\otimes 1. \notag
\end{gather}

\begin{theorem}[Theorem \ref{MAIN1}]\label{main2a}
The form $\varphi_{q,\ell}$ is closed. More precisely,
\[
d_{\calF}' \varphi_{q,\ell} =  d_{\calF}'' \varphi_{q,\ell} =  0
\]
and
\[
d_V \varphi_{q,\ell} =0. 
\]
\end{theorem}

\begin{proof}
First note that $d_{\calF}' \varphi_{q,\ell}=0$ is obvious from \eqref{dF}. 
From the `scalar- valued' case, see \cite{KMI}, we have $d_{\calF}
\varphi_{q,0} =d_{\calF}''\varphi_{q,0} =0$. In fact, this can be seen directly 
by \eqref{dF}, since one easily checks
\begin{equation}\label{dvan}
 \sum_{\alpha} (z_{\alpha} \otimes A(\omega_{\alpha\mu}) \otimes 1) \, \varphi_
{q,0} =0
\end{equation}
for any $\mu$. We then easily see
\[
d_{\calF}' \varphi_{q,\ell} = 
\left( d_{\calF}' \varphi_{q,0} \right) \cdot \varphi_{0,\ell} =0.
\]
For the action of $d_V$, we first note
\begin{equation}\label{rhoaction}
(1 \otimes 1 \otimes \rho(X_{\alpha\mu}) ) \,\varphi_{0,\ell} = \frac{-i}{4\pi} 
\sum_{k=1}^{\ell} (z_{\alpha} \otimes 1 \otimes A_k(e_{\mu})) \, \varphi_
{0,\ell-1}.
\end{equation}
\eqref{d1} then implies
\[
d_V \varphi_{q,\ell} = \frac{-i}{4\pi} \sum_{\alpha,\mu}\sum_{k=1}^{\ell} (z_
{\alpha} \otimes A(\omega_{\alpha\mu}) \otimes 1) \, \varphi_{q,0} \cdot (1 
\otimes 1 \otimes  A_k(e_{\mu})) \, \varphi_{0,\ell-1} =0
\]
by \eqref{dvan}.
\end{proof}

\subsection{Recursion}

We will now show Theorem~\ref{MAIN2}(i), the recursive formula for the
cohomology class $[\varphi_{q,\ell}]$.

For $j \geq 1$, we define operators $A_j(\sigma)$ by
\[
A_j(\sigma) = i \sum_{\alpha} \left(  \frac{\partial}{\partial z_{\alpha}}  - 
\frac1{4\pi}z_{\alpha} \right) \otimes 1 \otimes A_j(e_{\alpha}) \, - \, 
i \sum_{\mu} \left(  \frac{\partial}{\partial z_{\mu}}  - \frac1{4\pi}z_{\mu} 
\right) \otimes 1 \otimes A_j(e_{\mu}).
\]
We write $A(\sigma) = A_1(\sigma)$, and by \eqref{inter1} we note that this is 
the image in the Fock model of the operator $A(\sigma)$ in the Schroedinger 
model.



For $j \geq 1$, we introduce operators $h'_j$ by
\[
h'_j = \sum_{\alpha,\mu} \frac{\partial}{\partial z_{\alpha}} \otimes A^{\ast}
(\omega_{\alpha\mu}) \otimes A_j(e_{\mu}).
\]
We write $h'=h'_1$. Here $A^{\ast}(\omega_{\alpha\mu})$ denotes the (interior)
multiplication with $X_{\alpha\mu}$, i.e.,
$A^{\ast}(\omega_{\alpha\mu}) (\omega_{\alpha'\mu'}) =
\delta_{\alpha\alpha'}\delta_{\mu\mu'}$. 
We define a $(q-1)$-form $\Lambda^{(j)}_{q,\ell}$ by 
\[
\Lambda^{(j)}_{q,\ell} = \frac{-i}{p+q+\ell-1}h'_j \varphi_{q,\ell}.
\]
We write $\Lambda_{q,\ell} = \Lambda^{(1)}_{q,\ell}$.

\begin{theorem}[Theorem~\ref{MAIN2}]\label{main3a}
For any $1\leq j \leq \ell$, we have 
\[
\varphi_{q,\ell} = A_j(\sigma) \varphi_{q,\ell-1} \, + \, d \Lambda^{(j)}_
{q,\ell-1}  \, + \, \frac1{4\pi} \sum_{k=1}^{\ell-1} A_{jk}(f) \varphi_{q,\ell-
2}.
\]
\end{theorem}

For the proof of Theorem~\ref{main3a} we first compute $A_j(\sigma) \varphi_
{q,\ell-1}$:


\begin{lemma}\label{lem3a}
For any $1\leq j \leq \ell$, we have
\[
A_j(\sigma) \varphi_{q,\ell-1} =  \varphi_{q,\ell} 
\; + A_j \; + \; B_j  \;+ \; C_j^+. 
\]
Here
\begin{gather*}
A_j = \frac{i}{4\pi} \left( \sum_{\mu=p+1}^{p+q} z_\mu \otimes 1 \otimes A_j(e_
{\mu})\right) \varphi_{q,\ell-1}, \\
B_j= i  \sum_{\alpha=1}^p  (\frac{\partial}{\partial z_{\alpha}} \otimes 1 
\otimes 1) \,  \varphi_{q,0}  \cdot (1\otimes 1 \otimes A_j(e_{\alpha}) ) \, 
\varphi_{0,\ell-1}, \\
C_j^+ = \frac{1}{4\pi}\sum_{k=1}^{\ell-1} A_{jk}(f_+)\varphi_{q,\ell-2},
\end{gather*}
where $A_{jk}(f_+)$ is the insertion $\sum_{\alpha=1}^p A_j(e_{\alpha}) A_k(e_
{\alpha})$ in the $j$-th and $k$-th position in $T(V)$.
\end{lemma}

\begin{proof}
We write $A_j(\sigma) = A_j'(\sigma) + A_j''(\sigma)$ with $A_j'(\sigma) =  
\tfrac{-i}{4\pi} \sum_{\alpha} z_{\alpha} \otimes 1 \otimes A_j(e_{\alpha}) + 
\tfrac{i}{4\pi} \sum_{\mu} z_{\mu} \otimes 1 \otimes A_j(e_{\mu})
$.
We immediately see
\[
A_j'(\sigma) \varphi_{q,\ell-1} =  \varphi_{q,\ell} \; + \; A_j.
\]
On the other hand, we have 
\begin{align*}
 A_j''(\sigma) \varphi_{q,\ell-1} &= i  \sum_{\alpha=1}^p (
\frac{\partial}{\partial z_{\alpha}} \otimes 1 \otimes A_j(e_{\alpha}) ) 
(\varphi_{q,0} \cdot \varphi_{0,\ell-1}) \\
& = B_j + i \varphi_{q,0} \cdot  \left( \sum_{\alpha=1}^p (
\frac{\partial}{\partial z_{\alpha}} \otimes 1 \otimes A_j(e_{\alpha}) ) 
\right) \, \varphi_{0,\ell-1}.
\end{align*}
But the last term is equal to $C_j^+$. Indeed, a little calculation gives 
\begin{equation}\label{diffaction}
  (\frac{\partial}{\partial z_{\alpha}} \otimes 1 \otimes 1)  \, \varphi_
{0,\ell-1} =
\frac{-i}{4\pi} \sum_{k=1}^{\ell-1} ( 1 \otimes 1 \otimes A_k(e_{\alpha})) \, 
\varphi_{0,\ell-2},
 \end{equation}
from which the claim follows.
\end{proof}

Therefore Theorem~\ref{main3a} will follow from

\begin{proposition}\label{prop3a}
For any $1\leq j \leq \ell$, we have
\[
d \Lambda^{(j)}_{q,\ell-1} = -(A_j + B_j + C_j^-),
\]
where
\[
C_- = \frac{1}{4\pi}\sum_{k=1}^{\ell-1} A_{jk}(f_-)\varphi_{q,\ell-2}
\]
with $A_{jk}(f_-)=\sum_{\mu=1}^p A_j(e_{\mu}) A_k(e_{\mu})$.
\end{proposition}

\begin{proof}

Since $\varphi_{q,\ell-1}$ is closed, we have
\[
(p+q+\ell-2) d \Lambda_{q,\ell-1} = dh'_j \varphi_{q,\ell-1} = \{d,h_j'\} 
\varphi_{q,\ell-1},
\]
where $\{A,B\}$ denotes the anticommutator $AB + BA$. It is easy to see that
\[
\{d_{\mathcal{F}}',h_j'\} 
\varphi_{q,\ell-1} =0,
\]
so that we only need to compute 
$\{d_{\mathcal{F}}'',h_j'\} $ and $\{d_V,h_j'\} $.

\begin{lemma}\label{lemma3b}
As operators on 
 $\mathcal{F} \otimes \bigwedge{^{\ast}}(\mathfrak{p}^{\ast}) \otimes T(V)$,

\begin{equation}
4 \pi i \{d_{\mathcal{F}}'',h_j'\} = \sum_{\alpha,\mu} z_{\mu} \frac{\partial}
{\partial z_{\alpha}} z_{\alpha} \otimes 1 \otimes A_j(e_{\mu}) \; - \;
\sum_{\mu,\nu} z_{\mu} \otimes D_{\mu\nu} \otimes A_j(e_{\nu}). \tag{i}
\end{equation}

\begin{equation}
i \{d_V,h_j'\} = \sum_{\alpha,\mu} \frac{\partial}{\partial z_{\alpha}} \otimes 
1 \otimes A_j(e_{\mu}) \rho(X_{\alpha\mu}) \; + \; \sum_{\alpha,\beta}  \frac
{\partial}{\partial z_{\beta}} \otimes D_{\alpha\beta} \otimes A_j(e_{\alpha}). 
\tag{ii}
\end{equation}
Here the operators $ D_{\alpha\beta}$ and  $D_{\mu\nu}$ are the derivations of 
$\bigwedge{^{\ast}}(\mathfrak{p}^{\ast})$ determined by
\[
D_{\alpha\beta} \omega_{\gamma\mu} = \delta_{\beta\gamma} \omega_{\alpha\mu} 
\qquad \text{and} \qquad  D_{\mu\nu} \omega_{\alpha\lambda} = \delta_
{\lambda\nu} \omega_{\alpha\mu}.
\]

\end{lemma}

\begin{proof}



For (i), using the definitions of the operators and $\frac{\partial}{\partial z_
{\beta}} z_{\alpha} = z_{\alpha}\frac{\partial}{\partial z_{\beta}} + \delta_
{\alpha\beta}$, we easily see
\[
4\pi i \{d_{\mathcal{F}}'',h_j'\} = \sum_{\substack{\alpha,\mu \\ \beta,\nu}} z_
{\mu}z_{\alpha} \frac{\partial}{\partial z_{\beta}} \otimes \{ A_{\alpha\mu},A^
{\ast}_{\beta\nu}\} \otimes A_j(e_{\nu}) \;+\; \delta_{\alpha\beta} z_{\mu} 
\otimes A^{\ast}_{\beta\nu} A_{\alpha\mu} \otimes A_j(e_v).
\]
Here and in the following we write $A_{\alpha,\mu}$ for $A(\omega_{\alpha,\mu})
$. The Clifford identities imply
\[
\{ A_{\alpha\mu},A^{\ast}_{\beta\nu}\} = \delta_{\alpha\beta} \delta_{\mu\nu}
\]
and
\[
 A^{\ast}_{\alpha\nu} A_{\alpha\mu} = 
\begin{cases}
I - A_{\alpha\mu} A^{\ast}_{\alpha\nu}, &\text{ if $\mu=\nu$} \\
-  A_{\alpha\mu} A^{\ast}_{\alpha\nu}, &\text{ if $\mu\ne\nu$}.
\end{cases}
\]
Note 
\begin{equation}\label{Cliff}
\sum_{\alpha}  A_{\alpha\mu} A_{\alpha\nu}^{\ast} = D_{\mu\nu} \qquad \text
{and} \qquad
\sum_{\mu} A_{\alpha\mu}A^{\ast}_{\beta\mu} = D_{\alpha\beta}.
\end{equation}
We therefore obtain
\[
\sum_{\alpha,\mu} z_{\mu}\left( z_{\alpha} \frac{\partial}{\partial z_{\alpha}}
+1\right) \otimes 1 \otimes A_j(e_{\mu}) \;-\; \sum_{\alpha,\mu,\nu} z_{\mu} 
\otimes A_{\alpha\mu} A_{\alpha\nu}^{\ast} \otimes A_j(e_{\nu}),
\]
and the assertion follows from \eqref{Cliff}.

For (ii), first note $\rho(X_{\alpha\mu}) A_j(e_{\nu}) = \delta_{\mu\nu} A_j(e_
{\alpha})+ A_j(e_v) \rho(X_{\alpha\mu})$. One then obtains 
\begin{align*}
\{d_V,h_j'\} &= \sum_{\substack{\alpha,\mu \\\beta,\nu}} 
\frac{\partial}{\partial z_{\beta}}\otimes  \{ A_{\alpha\mu},A^{\ast}_{\beta\nu}
\} \otimes A_j(e_v) \rho(X_{\alpha\mu}) \;+\; 
\frac{\partial}{\partial z_{\beta}}\otimes 
 A_{\alpha\mu}A^{\ast}_{\beta\nu} \otimes \delta_{\mu\nu} A_j(e_{\alpha}) \\
&= \sum_{\alpha,\mu} \frac{\partial}{\partial z_{\alpha}} \otimes 1 \otimes A_j
(e_{\mu}) \rho(X_{\alpha\mu}) \;+\; \sum_{\alpha,\beta} 
\frac{\partial}{\partial z_{\beta}}\otimes D_{\alpha\beta} \otimes A_j(e_
{\alpha}),
\end{align*}
by \eqref{Cliff}.
\end{proof}

Proposition~\ref{prop3a} now follows from

\begin{lemma}\label{lemma3c}

\begin{equation}
 \{d_{\mathcal{F}}'',h_j'\} \varphi_{q,\ell-1} = -(p+q+\ell-2) A_j. \tag{i}
\end{equation}
\begin{equation}
\{d_V,h_j'\} \varphi_{q,\ell-1} = -(p+q+\ell-2)(B_j +C_j^-). \tag{ii}
\end{equation}

\end{lemma}

\begin{proof}

We use Lemma~\ref{lemma3b}. For (i), we first have
\begin{align*}
\sum_{\alpha,\mu} \left(z_{\mu}\frac{\partial}{\partial z_{\alpha}} z_{\alpha}
\otimes 1 \otimes A_j(e_{\mu}) \right) \varphi_{q,\ell-1} &= (p+q+\ell-1) \left
(\sum_{\mu} z_{\mu} \otimes 1 \otimes A_j(e_{\mu}) \right) \varphi_{q,\ell-1} \\
& = -4\pi i (p+q+\ell-1) A_j,
\end{align*}
which follows immediately from
\begin{equation}\label{K'3}
\sum_{\alpha} \left( \frac{\partial}{\partial z_{\alpha}} z_{\alpha} \otimes 1 
\otimes 1 \right) \varphi_{q,\ell-1} = 
\sum_{\alpha} \left(  (z_{\alpha} \frac{\partial}{\partial z_{\alpha}} +1) 
\otimes 1 \otimes 1 \right) \varphi_{q,\ell-1} = (p+q+\ell-1)\varphi_{q,\ell-1}
\end{equation}
by \eqref{K'1}. Furthermore, by \cite{KM90}, Lemma~8.2 we have
\[
(1\otimes D_{\mu\nu} \otimes 1) \varphi_{q,0} = \delta_{\mu\nu} \varphi_{q,0},
\]
and therefore
\begin{align*}
\sum_{\mu,\nu} \left(z_{\mu} \otimes D_{\mu\nu} \otimes A_j(e_{\nu})  \right) 
\varphi_{q,\ell-1} &= 
\sum_{\mu,\nu} (1\otimes D_{\mu\nu} \otimes 1) \varphi_{q,0} \cdot 
  \left(z_{\mu} \otimes 1 \otimes A_j(e_{\nu})  \right) \varphi_{0,\ell}
\\
& = \left(\sum_{\mu} z_{\mu} \otimes 1 \otimes A_j(e_{\mu}) \right) \varphi_
{q,\ell-1} = -4\pi i A_j.
\end{align*}
 Lemma~\ref{lemma3b} now gives (i). 
For (ii), by \eqref{rhoaction} we first note
\begin{equation}
(1 \otimes 1 \otimes \rho(X_{\alpha\mu}) \, \varphi_{q,\ell-1} = \frac{-i}{4
\pi} \sum_{k=1}^{\ell-1} (z_{\alpha} \otimes 1 \otimes A_k(e_{\mu}) ) \varphi_
{q,\ell-2}.
\end{equation}
Thus, using \eqref{K'3}, we see
\begin{align*}
 \sum_{\alpha,\mu}  \left( \frac{\partial}{\partial z_{\alpha}} \otimes 1 
\otimes A_j(e_{\mu}) \rho(X_{\alpha\mu}) \right) \varphi_{q,\ell-1}
& = \frac{-i}{4\pi} \left( \sum_{k=1}^{\ell-1} A_{jk}(f_-) \right) (p+q+\ell-
2)  \varphi_{q,\ell-2} \\
& = -i  (p+q+\ell-2) C_j^-.
\end{align*}

Finally, by \cite{KM90}, Lemma~8.2, we have
\[
(1 \otimes  D_{\alpha\beta} \otimes 1) \varphi_{q,0} = (z_{\beta} 
 \frac{\partial}{\partial z_{\alpha}} \otimes 1 \otimes 1) \varphi_{q,0}.
\]
Hence
\begin{align}\label{formel}
\sum_{\alpha, \beta}   & \left( \frac{\partial}{\partial z_{\beta}} 
 \otimes  D_{\alpha\beta}  \otimes A_j(e_{\alpha}) \right) \varphi_{q,\ell-1} 
\\ &= \sum_{\alpha, \beta}  \left( \frac{\partial}{\partial z_{\beta}} z_
{\beta} \otimes 1 \otimes  A_j(e_{\alpha}) \right) (  (\frac{\partial}{\partial 
z_{\alpha}}
 \otimes 1 \otimes 1) \varphi_{q,0} ) \cdot \varphi_{0,\ell-1} \notag \\ \notag
&= \sum_{\alpha, \beta}
 \left( \left( \frac{\partial}{\partial z_{\beta}} z_{\beta} \frac{\partial}
{\partial z_{\alpha}} \otimes 1 \otimes  1 \right) \, \varphi_{q,0} \right) 
\cdot ( 1 \otimes 1 \otimes A_j(e_{\alpha})) \, \varphi_{0,\ell-1} \\ \notag
& \quad +  \sum_{\alpha, \beta} 
( (\frac{\partial}{\partial z_{\alpha}}
 \otimes 1 \otimes 1) \varphi_{q,0} )\cdot 
\left( \frac{\partial}{\partial z_{\beta}} z_{\beta} \otimes 1 \otimes  A_j(e_
{\alpha}) \right) \varphi_{0,\ell-1},
\end{align}
The second term, using \eqref{K'3}, is equal to $-i(\ell-1) B$, while for the 
first we have
\begin{align*}
 \sum_{\alpha, \beta} \left( \left( \frac{\partial}{\partial z_{\beta}} (\frac
{\partial}{\partial z_{\alpha}} z_{\beta} - \delta_{\alpha\beta}) \otimes 1 
\otimes  1 \right) \, \varphi_{q,0} \right) \cdot ( 1 \otimes 1 \otimes A_j(e_
{\alpha})) \, \varphi_{0,\ell-1},
\end{align*}
and this, using \eqref{K'3} again, equals to $-i(p+q-1) B$. This finishes the 
proof of (ii).
\end{proof}

This concludes the proof of Proposition~\ref{prop3a} and hence the proof of 
Theorem~\ref{main3a}!
\end{proof}

\subsection{Holomorphicity}

We will now show Theorem~\ref{MAIN4}(i), i.e., that the cohomology class
 $[\varphi_{q,[\ell]}]$ is holomorphic.

Following \cite{KM90}, we define another operator $h$ on $\mathcal{F} \otimes 
\bigwedge ^{\ast}(\mathfrak{p}^{\ast}) \otimes S^{\ast}(V)$  by
\begin{gather*}
h =  \sum_{\alpha,\mu} z_{\mu} \frac{\partial}{\partial z_{\alpha}}
\otimes A^{\ast}(\omega_{\alpha\mu}) \otimes 1. 
\end{gather*}

\begin{definition}
We introduce a $(q-1)$-form $\psi_{q,\ell}$ by
\[
\psi_{q,\ell} := \frac{-1}{2(p+q -1)} \left(h \varphi_{q,0}\right) \cdot 
\varphi_{0,\ell}.
\]
\end{definition}
It will be convenient to note that we could have also defined $\psi_{q,\ell}$ 
by letting $h$ act on $\varphi_{q,\ell}$.
\begin{lemma}\label{newdef}
\[
\psi_{q,\ell} = \frac{-1}{2(p+q +\ell-1)} h \varphi_{q,\ell}.
\]
\end{lemma}

\begin{proof}
We have
\begin{align}\label{calc}
h \varphi_{q,\ell} & =  \left(h \varphi_{q,0}\right) \cdot \varphi_{0,\ell} 
+ 
\sum_{\alpha,\mu} (z_{\mu}  \otimes A^{\ast}(\omega_{\alpha\mu}) \otimes 1) 
\varphi_{q,0} 
\cdot (\frac{\partial}{\partial z_{\alpha}}
\otimes 1 \otimes 1) \varphi_{0,\ell} \\
&=  \left(h \varphi_{q,0}\right) \cdot \varphi_{0,\ell} - \frac{i}{4\pi} \sum_
{j=1}^{\ell} h''_j \varphi_{q,\ell-1},
\end{align}
with $h''_j = \sum_{\alpha,\mu} z_{\mu} \otimes A^{\ast}(\omega_{\alpha\mu}) 
\otimes A_j(e_{\alpha})$. But now one easily checks that
\[
\frac{-1}{4\pi}  h''_j \varphi_{q,\ell-1} = \varphi'_{q,0} \cdot \varphi_
{0,\ell}\]
with 
\[
\varphi'_{q,0} = \sum_{\mu, \alpha_1,\dots,\alpha_{q-1}} (-1)^{\mu-p-1} z_{\mu} 
z_{\underline{\alpha}} \otimes \omega_{\alpha_1 p+1} \wedge \cdots \wedge 
\widehat{\omega_{\alpha_{\mu-p} \mu}} \wedge \cdots  \omega_{\alpha_{q-1} p+q} 
\otimes 1.
\]
On the other hand, we also compute
\[
h \varphi_{q,0} = (p+q-1) \varphi'_{q,0}.
\]
>From this the lemma follows.
\end{proof}

\begin{theorem}[Theorem~\ref{MAIN4}]\label{mt4}
The action of the lowering operator $L$ on  $\varphi_{q,\ell}$ is given by 
\[
\omega(L) \varphi_{q,\ell} = d\left(\psi_{q,\ell} + \frac{1}2 \sum_{j=1}^{\ell} 
\Lambda^{(j)}_{q,\ell-1} \right) - \frac{1}{4 \pi} A(f)\varphi_{q,\ell-2}.
\]
\end{theorem}

\begin{remark}
This theorem is the generalization of one of the main points in \cite{KM90} for 
$\ell=0$, the trivial coefficient case. Namely, \cite{KM90}, Lemma~8.3 states
\begin{equation}\label{psiformel}
\omega(L) \varphi_{q,0} = d \psi_{q,0}.
\end{equation}
\end{remark}

\begin{proof}[Proof of Theorem~\ref{mt4}]

We first compute the left hand side:

\begin{lemma}\label{lemma4a}
\[
\omega(L) \varphi_{q,\ell} = \left(\omega(L)\varphi_{q,0}\right) \cdot \varphi_
{0,\ell} \; - \;  \sum_{j=1}^\ell B_j \; -  \;  \frac{1}{4 \pi} A(f_+)\varphi_
{q,\ell-2}, 
\]
with $B_j$ as in Lemma~\ref{lem3a}.

\end{lemma}

\begin{proof}
By Lemma~\ref{Fock2} we have
\begin{align}
\omega(L) \varphi_{q,\ell} &= \frac{-1}{8\pi} \sum_{\mu} \left(z_{\mu}^2 
\otimes 1 \otimes 1\right)  \varphi_{q,\ell}  \label{t1} \\
& \quad + 2\pi c_{q,\ell} \sum_{\substack{\g \\ \underline{\alpha}, \underline
{\beta}}} \left(\frac{\partial^2}{\partial z_{\g}^2} z_{\underline{\alpha}} 
\right) z_{\underline{\beta}} \otimes \omega_{\underline{\alpha}} \otimes e_
{\underline{\beta}}\label{t2}\\
& \quad + 4 \pi c_{q,\ell} \sum_{\substack{\g \\ \underline{\alpha}, \underline
{\beta}}} \left(\frac{\partial}{\partial z_{\g}} z_{\underline{\alpha}} \right) 
\left(
\frac{\partial}{\partial z_{\g}}z_{\underline{\beta}} \right) \otimes \omega_
{\underline{\alpha}} \otimes e_{\underline{\beta}} \label{t3} \\
& \quad + 2\pi c_{q,\ell} \sum_{\substack{\g \\ \underline{\alpha}, \underline
{\beta}}}  z_{\underline{\alpha}} \left(\frac{\partial^2}{\partial z_{\g}^2} z_
{\underline{\beta}} \right)\otimes \omega_{\underline{\alpha}} \otimes e_
{\underline{\beta}}. \label{t4}
\end{align}
The first two terms ((\ref{t1}) and (\ref{t2}) give $\left(\omega(L)\varphi_
{q,0}\right) \cdot \varphi_{0,\ell}$. By \eqref{diffaction}, the third term 
(\ref{t3}) is equal to
\[
-i \sum_{\g} \left( \frac{\partial}{\partial z_{\g}} \otimes 1 \otimes 1 
\right) \, \varphi_{q,0} \cdot \sum_{j=1}^{\ell} (1 \otimes 1 \otimes A_j(e_
{\g})) \varphi_{0,\ell-1} = - \sum_{j=1}^{\ell} B_j.
\]
For the  fourth term (\ref{t4}) in the sum above, we apply \eqref{diffaction} 
twice and obtain
\begin{multline*}
\frac{-i}{2} \varphi_{q,0} \cdot \sum_{j=1}^{\ell} \sum_{\g}\left( \frac
{\partial}{\partial z_{\g}} \otimes 1 \otimes A_j(e_{\g}) \right) \varphi_
{0,\ell-1} \\ = 
\frac{-1}{8\pi} \varphi_{q,0} \cdot \sum_{j=1}^{\ell} \sum_{k=1}^{\ell-1} \sum_
{\g}(1 \otimes 1\otimes  A_j(e_{\g}) A_k(e_{\g})) \varphi_{0,\ell-2} = \frac{-1}
{4\pi} A(f_+)  \varphi_{0,\ell-2}.
\end{multline*}
\end{proof}

We now compute $d \psi_{q,\ell}$:

\begin{lemma}\label{lemma4b}
\begin{equation}
d_{\mathcal{F}} \psi_{q,\ell} = \left(d_{\mathcal{F}} \psi_{q,0}\right)\cdot 
\varphi_{0,\ell} - \frac{1}2 \sum_{j=1}^{\ell}B_j. \tag{i}
\end{equation}
\begin{equation}
d_{V} \psi_{q,\ell} =   \frac{1}2 \sum_{j=1}^{\ell}A_j. \tag{ii}
\end{equation}
with $A_j$ and $B_j$ as in Lemma~\ref{lem3a}.

\end{lemma}

\begin{proof}

For (i), we first observe
\begin{equation}\label{l4b1}
d''_{\mathcal{F}} \left(  \left(h \varphi_{q,0}\right) \cdot \varphi_{0,\ell} 
\right)
= \left( d''_{\mathcal{F}}h \varphi_{q,0}\right) \cdot \varphi_{0,\ell}.
\end{equation}
For $d'_{\mathcal{F}}$ we get
\begin{align}
 d'_{\mathcal{F}}  \left(  \left(h \varphi_{q,0}\right) \cdot \varphi_{0,\ell} 
\right) &= -4\pi \sum_{\alpha, \mu} \left( \frac{\partial}{\partial z_{\alpha}} 
\otimes 1 \otimes 1 \right) \left[ \left( \left(
 \frac{\partial}{\partial z_{\mu}} \otimes A_{\alpha\mu} \otimes 1 \right)  h 
\varphi_{q,0} \right) \cdot \varphi_{0,\ell} \right] \notag \\
&= 
\left( d'_{\mathcal{F}} h \varphi_{q,0} \right) \cdot  \varphi_{0,\ell} \label
{l4b2}  \\
& \quad - 
4 \pi \sum_{\alpha, \mu} 
 \left( \left(
 \frac{\partial}{\partial z_{\mu}} \otimes A_{\alpha\mu} \otimes 1 \right)  h 
\varphi_{q,0} \right) \cdot \left(  \left(\frac{\partial}{\partial z_{\alpha}} 
\otimes 1 \otimes 1 \right) \varphi_{0,\ell} \right). \label{l4b3}
\end{align}
For the first term in (\ref{l4b3}) we see
\begin{align*}
\sum_{ \mu} 
  \left(
 \frac{\partial}{\partial z_{\mu}} \otimes A_{\alpha\mu} \otimes 1 \right) \left
(h \varphi_{q,0} \right) &= 
\sum_{\substack{\mu,\nu \\ \beta}} \left( \frac{\partial}{\partial z_{\beta}}
\frac{\partial}{\partial z_{\mu}} z_{\nu}\otimes A_{\alpha\mu} A_{\beta \mu}^
{\ast} \otimes 1 \right) \varphi_{q,0} \\
& = \sum_{\beta, \mu} \left( \frac{\partial}{\partial z_{\beta}} \otimes  A_
{\alpha\mu} A_{\beta \nu}^{\ast} \otimes 1 \right) \varphi_{q,0} \\
& = \sum_{\beta} \left( \frac{\partial}{\partial z_{\beta}} \otimes D_
{\alpha\beta} \otimes 1 \right) \varphi_{q,0}.
\end{align*}
Combining this with \eqref{diffaction} we obtain for (\ref{l4b3})
\[
i \sum_{j=1}^{\ell} \sum_{\alpha,\beta} \left(\left( \frac{\partial}{\partial z_
{\beta}} \otimes D_{\alpha\beta} \otimes A_j(e_{\alpha}) \right) \varphi_{q,0} 
\right)\cdot \varphi_{0,\ell-1}.
\]
But this is exactly (up to a constant) the term (\ref{formel}), i.e., 
(\ref{l4b3}) is equal to $ (p+q-1)\sum_{j=1}^{\ell} B$. This together with (\ref
{l4b1}),(\ref{l4b2}) and collecting the constants implies (i). 

For (ii), we will use Lemma~\ref{newdef}. We easily see
\[
\{d_V,h\} = \sum_{\alpha,\mu} z_{\mu} \frac{\partial}{\partial z_{\alpha}} 
\otimes 1 \otimes \rho(X_{\alpha\mu}).
\]
By \eqref{rhoaction} and \eqref{K'1}, we get 
\begin{align*}
d_V h \varphi_{q,\ell} &= \{d_V,h\} \varphi_{q,\ell} \\
&=  \frac{-i}{4\pi} \sum_{j=1}^{\ell} 
 \sum_{\alpha,\mu} \left( z_{\mu} \frac{\partial}{\partial z_{\alpha}} z_
{\alpha} \otimes 1 \otimes A_j(e_{\mu}) \right) \varphi_{q,\ell-1} \\
&=  -   \sum_{j=1}^{\ell} (p+q+\ell-1) A_j.
\end{align*}
This implies (ii).
\end{proof}

We are now in the position to finish the proof of Theorem~\ref{mt4}: Combining 
Lemma~\ref{lemma4a}, Lemma~\ref{lemma4b} and Proposition~\ref{prop3a}, we get
\begin{align*}
 \omega(L) \varphi_{q,\ell} -  d\left(\psi_{q,\ell} + \frac{1}2 \sum_{j=1}^q 
\Lambda^{(j)}_{q,\ell-1} \right)   & = 
(\omega(L) \varphi_{q,0}) \cdot \varphi_{0,\ell} - \sum_{j=1}^{\ell} B_j - \frac
{1}{4\pi} A(f_+)\varphi_{q,\ell-2} \\ &\quad - 
(d_{\mathcal{F}} \psi_{q,0}) \cdot \varphi_{0,\ell} + \frac{1}2\sum_{j=1}^
{\ell} B_j  - \frac{1}2 \sum_{j=1}^{\ell} A_j \\ & \quad + \frac{1}2 \sum_{j=1}^
{\ell}(A_j +B_j +C_j^-) \\
& = - \frac{1}{4 \pi} A(f) \varphi_{q,\ell-2},
\end{align*}
via (\ref{psiformel}) and $ A(f_-) \varphi_{q,\ell-1}= \frac{1}2 \sum_{j=1}^
{\ell} C_j^-$, since $A(f) = A(f_+) - A(f_-)$.
\end{proof}

\section{Main Result}\label{Theta}

In this section, we first construct the cohomology class $[\theta_{nq,[\la]}]$ 
and use the fundamental properties of $\varphi_{nq,[\la]}$ to derive our main 
result.

First note that over $\R$, the Weil representation action of 
 the standard Siegel parabolic in $Sp(n,\R)$  on $\varphi \in \calS(V^n)$ is 
given by 
\begin{equation*}
\omega \left( \left( \begin{smallmatrix} a&0\\0&^ta^{-1}
    \end{smallmatrix} \right) \right) \varphi(\x) =
 (\det a)^{m/2}
  \varphi(\x a)
 \end{equation*}
 for $a \in GL^+_n(\R)$, and
\begin{equation*}
\omega \left( \left( \begin{smallmatrix} 1&b\\0&1
    \end{smallmatrix} \right) \right) \varphi(\x) =
e^{\pi i tr(b(\x,\x))} \varphi(\x)
\end{equation*}
for $b \in Sym_n(\R)$. Here $\x =(x_1, \cdots, x_n) \in V^n$, as before.

Globally, we let $\mathbb{A} = \mathbb{A}_{\mathbb{K}}$ be the ring of adeles 
of $\mathbb{K}$. Let $ G'(\mathbb{A})$ be the two-fold cover $Sp(n,\mathbb{A})
$, which acts on $S(\underline{V}^n(\mathbb{A}))$ via the (global) Weil 
representation  $\omega = \omega_{\underline{V}}$.

For $g' \in G'(\mathbb{A})$, we define the standard theta kernel associated to 
a Schwartz function $\varphi \in S(V^n_{\A})$ by 
\[ 
\theta(g',\varphi):=  \sum_{\x \in \underline{V}^n} \omega(g')\varphi(\x).
\]
By abuse of notation, we let $\varphi_{nq,[\la]} = \otimes_{i=1}^r \varphi_{v_i}
$, where $\varphi_{v_i }$ is the Schwartz form $\varphi_{nq,[\la]}$ at the 
first infinite place and the Standard Gaussian $\varphi_{0}$ at the other 
places, and we define $\varphi_{nq,\ell}$ in the same way. For the finite 
places, we let $\varphi_f$ correspond to the characteristic function of $h + 
\mathfrak{b}L^n$.

Given $\tau =( \tau_1,\dots,\tau_r) \in \h_n^r$ we $g'_{\tau} \in
Sp(n,\mathbb{K}_{\infty})$  be a standard element  which moves the base point 
$(i,\dots,i) \in \h_n^r$ to $\tau$, i.e.,
\[
g'_{\tau} = \begin{pmatrix}  1&u \\ 0&1 \end{pmatrix} \begin{pmatrix} a&0 \\ 0& 
{^ta^{-1}} \end{pmatrix},
\]
with $v = a {^ta^{-1}}$. We consider $g'_{\tau} \in G'(\mathbb{A})$ in the 
natural way.

We write $\rho_{\la}$ for the representation action  of $GL_n(\C)$ on $S_{\la}
(\C^n)$.

For $\tau \in \h_n^r$, we then define
\begin{align*}
\theta_{nq,[\la]}(\tau,z)  &=     N^{\mathbb{K}}_{\Q}(\det(a)^{-m}) \rho^{\ast}_
{\la}(a) \theta(g'_{\tau},\varphi_{\infty}) \\
&=\sum_{\x \in h + \mathfrak{b}L^n} \rho_{\la}^{\ast}(a) \varphi_{nq,[\la]}(\x 
a) e_{\ast} ((\x,\x)u/2), 
\end{align*}
and similarly, $\theta_{nq,\ell}(\tau,z)$. Here 
\[
e_{\ast}(A) = \exp\left(2\pi i \sum_{j=1}^r tr(\la_i(A)) \right),
\]
 for $A \in M_{n,n}(\mathbb{K})$. 

(Slightly abusing e denote by $\mathcal{A}(\G',S_{\la}(\C^n)^{\ast} \otimes 
\det^{-m/2})$ the space of vector-valued (not necessarily holomorphic) Hilbert-
Siegel modular forms of degree $n$ for a congruence subgroup $\G'$ and for the 
representation $( \rho_{\la}^{\ast} \otimes\det^{-m/2}, \det^{-m/2}, 
\cdots,\det^{-m/2})$.

Then Theorems~\ref{MAIN3} and the standard theta machinery give us

\begin{proposition}
 \[
\theta_{nq,[\la]}(\tau,z) \in \mathcal{A}(\G',S_{\la}(\C^n) \otimes \det{^
{m/2}}) \otimes A^{nq}(M,S_{[\la]} \calV),
\]
i.e, $\theta_{nq,[\la]}(\tau,z)$ is a vector-valued non-holomorphic Hilbert-
Siegel modular form for the  representation $( \rho_{\la}^{\ast}\otimes \det^{-
m/2}, \det^{-m/2}, \dots,\det^{-m/2})$ with values in the ${ S_{[\la]}(\calV)  }
$-valued closed differential $nq$-forms of the manifold $M$. 

\end{proposition}

The Fourier expansion of  $\theta_{nq,[\la]}(\tau)$ is given by
\[
\theta_{nq,[\la]}(\tau) = \sum_{ \beta \in
  Sym_n(\mathbb{K})}  \theta_{\beta, nq,[\la] }(v)
e_{\ast}( \beta \tau),
\]
with
\[
 \theta_{\beta,  nq,[\la] }(v) = \sum_{\x \in \mathcal{L}_{\beta}}
\rho_{\la}^{\ast}(a)\varphi_{nq,[\la]}(\x a) e_{\ast}((\x,\x)v/2 ).
\]

\begin{definition}
If $\eta$ is a rapidly decreasing ${ S_{[\la]}({\calV})  }$-valued closed $(p-n)
q$ form
on $M$ representing a class $[\eta] \in H_c^{(p-n)q}(M, S_{[\la]}
\calV)$, we define
\[
\Lambda_{nq,[\la]}(\tau, \eta) = \int_M \eta \wedge \theta_{nq,[\la]}(\tau) \in 
S_{\la}(\C^n)^{\ast}.
\]
We define $\Lambda_{nq,\ell}(\tau, \eta)$ in the same fashion for $\eta$ taking 
values in $T^{\ell}(\calV)$. If $\eta$ is  ${ S_{[\la]}({\calV})  }$-valued, 
then we have
\[
 \Lambda_{nq,\ell}(\tau, \eta) =\Lambda_{nq,[\la]}(\tau, \eta).
\]

\end{definition}

Before we can state our main result, we need a bit more
notation. For $q=2k$ even, we let $e_q$ be the Euler form of the symmetric 
space $D$ (which is the Euler class  of the tautological vector bundle over 
$D$, i.e., the fiber over a point $z \in D$ is given by the negative $q$-plane 
$z$) and zero for $q$ odd. Here $e_q$ is normalized such that it is given in 
$\bigwedge^q(\mathfrak{p}^{\ast})$ by
\begin{equation*}
e_q = \left(-\frac{1}{4\pi} \right)^k \frac{1}{k!} \sum_{\sigma \in S_q} \sgn
(\sigma) \Omega_{p+\sigma(1),p+\sigma(2)} \dots  \Omega_{p+\sigma(2k-1),p+\sigma
(2k)},
\end{equation*}
with
\begin{equation*}
\Omega_{\mu\nu} = \sum_{\alpha=1}^p \omega_{\alpha\mu} \wedge \omega_
{\alpha\nu}.
\end{equation*}

\begin{remark}\label{ThomB}
The main result of \cite{KM90} is that in the scalar valued case, the
generating series
\[
 \sum_{ t=0}^n
\sum_{\substack{\beta \geq 0 \\ \rank \beta = t} } \left( PD(C_{\beta})
\wedge e_q^{n-t} \right) \, e_{\ast}( \beta \tau),
\]
is a classical Hilbert-Siegel modular form of weight $m/2$. The key point is 
that for nondegenerate $\x$ such that $(\x,\x)$ positive semidefinite of rank 
$t$, $\varphi_{nq,0}(\x)$ is (essentially) a Thom form for the cycle $C_X \cap 
e_q^{n-t}$. One has
\[
\int_{\G_X \back D} \eta \wedge \varphi_{nq,0}(\x a) =\left(\int_{C_X} \eta 
\wedge e_q^{n-t}\right) e_{\ast}(-(\x,\x)v /2).
\]
\end{remark}

\begin{remark}
Actually, in \cite{KM90} the noncompact hyperbolic case of signature $(p,1)$ 
with $n=p-1$
is excluded (when the cycles are infinite geodesics). In the following, we will
also exclude this case. Note however that for signature $(2,1)$ and
$n=1$ this restriction was removed in \cite{FM1}, and our result will
also hold in that particular case.
\end{remark}

The following two theorems are the generalization of the main result of \cite
{KM90}.

\begin{theorem}\label{MainFinal1}
The cohomology class $[\theta_{nq,[\la]}]$ is holomorphic; i.e., it
defines a holomorphic Siegel modular form of genus $n$ with values in
$S_{\la}(\C^n)^{\ast} \otimes \chi(-m/2)$ and with coefficients in $  H^{nq}
(M,S_{[\la]} \calV)$. 
\end{theorem}

\begin{proof}
This follows immediately from Theorem~\ref{MAIN4}: We have $[\overline{\partial}
\varphi_{nq,[\la]}]=0$, thus [$\overline{\partial}\theta_{nq,[\la]}](\tau) =0$.
\end{proof}

\begin{theorem}\label{MainFinal}
The Fourier expansion of 
$[\theta_{nq,[\la]}](\tau)$ is given by
\[
[\theta_{nq,\la}](\tau) =  \sum_{ t=0}^n
\sum_{\substack{\beta \geq 0 \\ \rank \beta = t} } \left( PD(C_{\beta,[\la]})
\wedge e_q^{n-t} \right) \, e_{\ast}( \beta \tau),
\]
where $ PD(C_{\beta,[\la]})$ denotes the Poincar\'e dual class of
$PD(C_{\beta,[\la]})$. Furthermore, if $q$ is odd or $i(\la)=n$, then 
$[\theta_{nq,\la}](\tau)$ is a cusp form.

\end{theorem}

This is equivalent to

\begin{theorem}
For $\eta$ a rapidly decreasing closed $q(n-p)$ form on $M$ with values in $S_
{[\la]}(\calV)$, the generating series 
\[
\Lambda_{nq,[\la]}(\tau,\eta) = \sum_{ t=0}^n
\sum_{\substack{\beta \geq 0 \\ \rank \beta = t} }
\int_{C_{\beta,[\la]}} \left(\eta \wedge e_q^t \right) \,  e_{\ast}(
\beta \tau)
\]
is a holomorphic Siegel modular form of type
$(\det^{m/2} \otimes \rho_{\la}, \det^{m/2}, \cdots,\det^{m/2})$.

Note that by Theorem~\ref{labasis} a basis of $S_{\la}(\C^n)$ is given by 
$ \pi_{\la}\eps_{f(\la)}$, where $f(\la)$ runs through the semistandard
fillings $SS(\la,n)$. With respect to this basis, we define the
$\pi_{\la}f(\la)$ component $(\Lambda_{[\la]}(\tau,\eta))_{\pi_{\la}f(\la)}$ by 
\[
(\Lambda_{nq,[\la]}(\tau,\eta))_{\pi_{\la} f(\la)} =  \left(\int_M \eta \wedge 
\theta_{nq,[\la]}(\tau)\right) \left( \eps_{f(\la)} \right).
\]
Note here that the value of $\Lambda$ at $\eps_{f(\la)}$ and $\pi_{\la}f(\la)$ 
is the same. We then have
\begin{align*}
(\Lambda_{nq,[\la]}(\tau,\eta))_{\pi_{\la}f(\la)} & = \sum_{ t=0}^n
\sum_{\substack{\beta \geq 0 \\ \rank \beta = t} }
\int_{C_{\beta,[f(\la)]}} \left(\eta  \wedge e_q^t \right) \,  e_{\ast}(
\beta \tau) \\
 & = \sum_{ t=0}^n
\sum_{\substack{\x \in \calL^c_t \\ (\x,\x) \geq 0}} 
\int_{C_{\x}} \left((\eta, {\x_{f(\la)}})  \wedge e_q^t \right) \,  e_{\ast} 
((\x,\x)\tau/2),
\end{align*}
where 
\[
\calL^c_t = \{ \x \in h + \mathfrak{b}L^n: \rank (\x,\x) = t;  \text{$\x$ 
nondegenerate}\}
\]
\end{theorem}

\begin{proof}

We denote the $\beta$ Fourier coefficient of $(\Lambda_{nq,[\la]}(\tau,\eta))_{f
(\la)}$ by 
\[
a_{\beta} = \left(\int_{M} \eta \wedge\sum_{\x \in \mathcal{L}_{\beta}}
\rho_{\la}^{\ast}(a)\varphi_{nq,[\la]}(\x a) \right) ((a^{-1}\eps)_{f(\la)}) e_
{\ast}(\beta v).
\]

We first note

\begin{lemma}
Assume that $\beta$ \emph{not} positive semidefinite. Then
\[
a_{\beta}=0.
\]
\end{lemma}

\begin{proof}
For $n>1$, this follows from the Koecher principle, since 
$\Lambda_{[\la]}(\tau,\eta)$ is holomorphic. For $n=1$, so that $\beta <0$, the 
vanishing follows from the vanishing in the trivial coefficient case by an 
argument similar to the positive definite coefficient, see Lemma~\ref{nondegen} 
below.
\end{proof}

For $\beta$ positive semidefinite, we write
\[
a_{\beta}^c = \left(\int_{M} \eta \wedge 
\sum_{\x \in \mathcal{L}^c_{\beta}}     \rho_{\la}^{\ast}(a)\varphi_{nq,[\la]}
(\x a) \right) ((a^{-1}\eps)_{f(\la)}) e_{\ast}(\beta v)
\]
for the contribution of the closed orbits and $a_{\beta}^d = a_{\beta} - a_
{\beta}^c$ for the degenerate  part.

\begin{lemma}\label{nondegen}
Assume $\beta$ is positive semidefinite of rank $t$. Then 
\[
a_{\beta}^c =\int_{C_{\beta,[f(\la)]}} \left(\eta \wedge e_q^{n-t} \right) 
\]
\end{lemma}

\begin{proof}

By the usual unfolding argument, we obtain
\begin{align*}
 a_{\beta}^c e_{\ast}(- \beta v) & = \left( \int_{\G \back D} \eta \wedge \sum_
{x \in \mathcal{L}^c_{\beta}} \varphi_{nq,[\la]}(\x a) \right)   ((a^{-1}\eps)_
{f(\la)})  \\
 & = \left( \int_{\G \back D} \eta \wedge \sum_{x \in \G \back \mathcal{L}^c_
{\beta}} \sum_{\g \in \G_X \back \G} \g^{\ast}  \varphi_{nq,[\la]}(\x a) 
\right)   ((a^{-1}\eps)_{f(\la)})  \\
&=  \sum_{x \in \G \back \mathcal{L}^c_{\beta}} \left(\int_{ \G_x \back D} \eta 
\wedge \varphi_{nq,[\la]}(\x a) \right)   ((a^{-1}\eps)_{f(\la)}).
\end{align*}
But now by Theorem~\ref{MAIN2} and Lemma~\ref{auxlemma}, we have
\[
 [\varphi_{nq,[\la]}(\x a ) ((a^{-1}\eps)_{f(\la)})] = \left[(1 \otimes 1 
\otimes \x_{[f(\la)]}) \varphi_{nq,0}(\x a)\right].
\]
Thus
\begin{align*}
\left(\int_{ \G_x \back D} \eta \wedge\rho_{\la}^{\ast}(a) \varphi_{nq,[\la]}
(\x a) \right)   (\eps_{f(\la)}) & = \int_{ \G_x \back D} \eta \wedge (1 
\otimes 1 \otimes \x_{[f(\la)]}) \varphi_{nq,0} (\x a) \\
&= \int_{ \G_x \back D} (\eta, \x_{[f(\la)]}) \wedge  \varphi_{nq,0}(\x a) \\
&= \left(\int_{C_X}   (\eta , \x_{[f(\la)]}) \wedge e_q^t\right) e_{\ast}(- 
(\x,\x) v/2), 
\end{align*}
by Remark~\ref{ThomB}. This implies $a_{\beta}^c = \int_{C_{\beta,[f(\la)]}} 
\eta \wedge e_q^t$, as claimed.
\end{proof}

It remains to show

\begin{lemma}
Assume $\beta$ is positive semidefinite. Then
\[
a_{\beta}^d =0
\]

\end{lemma}

\begin{proof}
The recursion formula reduces this to the analogous statement for the singular 
coefficients in the scalar-valued case in the same way as in Lemma~\ref
{nondegen}. The Lemma then follows from the vanishing of those coefficients in 
the scalar-valued case, see \cite{KM90}, \S4. We leave the details to the 
reader.
\end{proof}

This concludes the proof of the theorem.
\end{proof}

\appendix

\section{The Fock model}

We briefly review the construction of the Fock model of the (infinitesimal) 
Weil representation of the symplectic Lie algebra $\mathfrak{sp}(W \otimes \C)
$, where $(W, \langle \, , \, \rangle)$ denotes a non-degenerate real 
symplectic space of dimension $2N$. We follow \cite{Adams,KM90}. We let $J_0$ 
be a positive definite complex structure on $W$, i.e., the bilinear form given 
by $\langle w_1,J_0w_2\rangle$ is positive definite. Let 
$e_1,\dots,e_N;f_1,\dots,f_N$ be a standard symplectic basis of $W$ so that 
$J_0 e_j=f_j$ and $J_0 f_j = -e_j$. We decompose
\[
W \otimes \C = W' \oplus W''
\]
into the $+i$ and $-i$ eigenspaces under $J_0$. Then $w_j' = e_j-if_j$ and 
$w_j'' = e_j +i f_j$, $j =1,\dots,N$ form a basis for $W'$ and $W''$ 
respectively.
We identify $Sym^2(W)$ with $\mathfrak{sp}(W)$ via
\[
(x \circ y)(z) = \langle x,z \rangle y + \langle y,z \rangle x.
\]

Given $\la \in \C$, we define the quantum algebra (see \cite{Howe}) $\mathcal{W}
_{\la}$ to be the tensor algebra $T(W \otimes \C)$ modulo the two sided ideal 
generated by the elements of the form $x \otimes y - y \otimes x - \la \langle 
x,y \rangle 1$. We let $p: T(W \otimes \C) \to \mathcal{W}_{\la}$  be the 
quotient map.
Since $T(W \otimes \C)$ is graded, we have a filtration $F^{\ast}$ on $\mathcal
{W}_{\la}$, and we easily that $[F^k \mathcal{W}_{\la},
F^{k'} \mathcal{W}_{\la}] \subset F^{k+k'-2} \mathcal{W}_{\la}$. Hence $F^2 
\mathcal{W}_{\la}$ is a Lie algebra. Furthermore, we have a split extension of 
Lie algebras
\[
0 \longrightarrow F^1 \mathcal{W}_{\la} \longrightarrow 
 F^2 \mathcal{W}_{\la} \longrightarrow \mathfrak{sp}(W \otimes \C) 
\longrightarrow 0.
\]
Here the second map is $p(x \otimes y) \mapsto \la (x \circ y) \in Sym^2(W 
\otimes \C) \simeq \mathfrak{sp}(W \otimes \C)$, while the splitting map $j: 
Sym^2(W \otimes \C) \to 
 F^2 \mathcal{W}_{\la}$ is given by
\[
j(x \circ y) = \frac{1}{2\la} \left( p(x)p(y) + p(y)p(x) \right).
\]
(Note the sign error in \cite{KM90}, p.~151). 

We let $\mathcal{J}$ be the left ideal in $\mathcal{W}_{\la}$ generated by 
$W'$. 
 The projection $p$ induces an isomorphism of the symmetric algebra $Sym^{\ast}
(W'')$ with  $\mathcal{W}_{\la}/ \mathcal{J}$. We denote  by $\rho_{\la}$ the 
action of $\mathcal{W}_{\la}$ on $\mathcal{W}_{\la}/ \mathcal{J} \simeq Sym^
{\ast}(W'') $ given by left multiplication. We now identify $Sym^{\ast}(W'')$ 
with the polynomial functions $\mathcal{P}(\C^N)= \C[z_1,\dots,z_n]$ on $W'$ 
via $z_j(w''_k) = \langle w'_j, w''_k \rangle = 2i \delta_{jk}$ and observe 
that then the action of $W \subset \mathcal{W}_{\la}$ on $\mathcal{P}(\C^N)$ is 
given by
\[
\rho_{\la}(w''_j) = z_j \qquad \text{and} \qquad \rho_{\la}(w'_j) = 2 i \la 
\frac{\partial}{\partial z_j}.\]
This determines the action of  $\mathcal{W}_{\la}$, and we obtain an action 
$\omega_{\la}= \rho_{\la} \circ j$ of $\mathfrak{sp}(W \otimes \C)$ on $\mathcal
{P}(\C^N)$. 
This is Fock model of the Weil representation with central character $\la$.

\smallskip

We now let $V$ be a real quadratic space of signature $(p,q)$ (for the moment, 
we change notation and denote the standard basis elements by $v_{\alpha}$ and 
$v_{\mu}$), and let $W$ be a real symplectic space over $\R$ of dimension $2n$ 
(with standard symplectic basis $e_j$ and $f_j$, $j =1,\dots,n$). We consider 
the symplectic space $\W = V \otimes W $ of dimension $2n(p+q)$, and note that 
$\mathbb{J} = \theta \otimes J$ defines a positive definite complex structure 
on $\W$. Here $\theta$ is the Cartan involution with respect to the above basis 
of $V$, while $J$ is the 
 positive define complex structure with respect to the above symplectic basis of 
$W$. Then the $+i$-eigenspace $\W'$ of $\mathbb{J}$ is spanned by the $v_
{\alpha} \otimes w'_j$ and $v_{\mu} \otimes w_j''$, while the $-i$ eigenspace 
$\W''$ is spanned by the $v_{\alpha} \otimes w_j''$ and $v_{\mu} \otimes w_j'$.

We naturally have $\mathfrak{o}(V) \times \mathfrak{sp}(W) \subset \mathfrak{sp}
(V \otimes W)$, and one easily checks that the inclusions $j_1: \mathfrak{o}(V) 
\simeq \bigwedge^2(V) \to  \mathfrak{sp}(V \otimes W) \simeq Sym^2(V \otimes W)
$ and $j_2: \mathfrak{sp}(W)  \to  \mathfrak{sp}(V \otimes W) \simeq Sym^2(V 
\otimes W)$ are given by
\begin{align*}
j_1\left( v_1 \wedge v_2 \right) & = \frac{1}{2i} \left[ \sum_{j=1}^n (v_1 
\otimes w_j') \circ (v_2 \otimes w_j'') - \sum_{j=1}^n (v_1 \otimes w_j'') 
\circ (v_2 \otimes w_j') \right]  \\
j_2\left( w_1 \circ w_2 \right)& = \sum_{\alpha=1}^p (v_{\alpha} \otimes w_1) 
\circ (v_{\alpha}\otimes w_2)  - \sum_{\mu =p+1}^{p+q} (v_{\mu} \otimes w_1) 
\circ (v_{\mu}\otimes w_2)
\end{align*}
with $v_1,v_2 \in V$ and $w_1,w_2 \in W$, see \cite{KM90}, Lemma~7.3.

We write $\calF= \mathcal{P}(\C^{n(p+q)})$ for the Fock model $\calF$ of the 
(infinitesimal) Weil representation of $\mathfrak{sp}(V \otimes W)$. We denote 
the
variables in $\mathcal{P}(\C^{n(p+q)})$ by $z_{\alpha j}$ corresponding to 
 $v_{\alpha} \otimes w_j''$ and $z_{\mu j}$  corresponding to  $v_{\mu} \otimes 
w_j'$. For $n=1$ we drop the subscript $j(=1)$. We have
\begin{align*}
\rho_{\la}( v_{\alpha} \otimes w'_j) &= 2i\la \frac{\partial}{\partial z_
{\alpha j}},  & \rho_{\la} ( v_{\alpha} \otimes w_j'') & = z_{\alpha j}, \\
 \rho_{\la}( v_{\mu } \otimes w''_j) &= 
2i\la \frac{\partial}{\partial z_{\mu j}},  &\rho_{\la} ( v_{\mu} \otimes w_j') 
& = z_{\mu j}.
\end{align*}

>From this we easily obtain the following formulas for the action of $\mathfrak
{o}(V) \times \mathfrak{sp}(W)$ in $\calF$:
(the formulas differ from the ones given in \cite{KM90} by a sign due to the 
sign error mentioned above).

\begin{lemma}\label{Fock1}
For the orthogonal group $\mathfrak{o}(V) = \mathfrak{k} \oplus \mathfrak{p}$, 
we write $X_{rs} = v_r \wedge v_s \in  \bigwedge^2(V) \simeq \mathfrak{o}(V) $. 
So $ \mathfrak{k}$ is spanned by $X_{\alpha\beta}$ and $X_{\mu\nu}$,
while  $\mathfrak{p}$ is spanned by the $X_{\alpha\mu}$. Then
\begin{align*}
\omega(X_{\alpha\beta}) & = -  \sum_{j=1}^n z_{\alpha j} \frac{\partial}
{\partial z_{\beta j}}   - z_{\beta j} \frac{\partial}{\partial z_{\alpha j}}, 
\\
\omega(X_{\mu\nu }) & = 
\sum_{j=1}^n z_{\mu j} \frac{\partial}{\partial z_{\nu j}}   - z_{\nu j} \frac
{\partial}{\partial z_{\mu j}}, \\
\omega(X_{\alpha\mu}) & = 2i\la  \sum_{j=1}^n \frac{\partial^2}{\partial
  z_{\alpha j} \partial z_{\mu j}} - \frac{1}{2i\la} \sum_{j=1}^n z_{\alpha j}z_
{\mu j}.
\end{align*}

\end{lemma}

\begin{lemma}\label{Fock2}

For the symplectic group, we note that in the decomposition $\mathfrak{sp}(W 
\otimes \C) = \mathfrak{k}' \oplus \mathfrak{p}^+ \oplus \mathfrak{p}^-$, 
$\mathfrak{k}' \simeq \mathfrak{gl}_n{\C}$ is spanned by the elements of the 
form $w'_j \circ w''_k$,  $\mathfrak{p}^+$ is spanned by $ w''_j \circ w''_k$ 
and $\mathfrak{p}^-$ is spanned by $w'_j \circ w'_k$ ($1 \leq j,k \leq n$). Then
\begin{align*}
\omega(w'_j \circ w''_k) &= 2i \left[ \sum_{\alpha=1}^p z_{\alpha k} \frac
{\partial}{\partial z_{\alpha j}} - \sum_{\mu = p+1}^{p+q} z_{\mu j} \frac
{\partial}{\partial z_{\mu k}} \right] + i(p-q) \delta_{jk}, \\
\omega(w''_j \circ w''_k) &=  \frac1{\la}  \sum_{\alpha=1}^pz_{\alpha j}  z_
{\alpha k} + 4 \la  \sum_{\mu = p+1}^{p+q} \frac{\partial^2}{\partial z_{\mu j} 
\partial z_{\mu k}}, \\
\omega(w'_j \circ w'_k) &= -4\la  \sum_{\alpha=1}^p \frac{\partial^2}{\partial 
z_{\alpha j} \partial z_{\alpha k}} - \frac1{\la} \sum_{\mu = p+1}^{p+q} z_{\mu 
j} z_{\mu k}.
\end{align*}

\end{lemma}

Note that for  $n=1$, we have $\mathfrak{sp}(W \otimes \C) \simeq 
\mathfrak{sl}_2(\C)$, and (for $\la = 2 \pi i$) the action of 
$L :=  \tfrac12 \left( \begin{smallmatrix}1 & -i \\ -i & -1
\end{smallmatrix} \right) = \tfrac{i}{4} w'_1 \circ w'_1 $ and 
$R
:=  \tfrac12 \left( \begin{smallmatrix}
1 & i \\ i & -1
\end{smallmatrix} \right) = \tfrac{-i}{4}  w''_1 \circ w''_1$  correspond to 
the classical Maass lowering and raising operators on the upper half plane.

\smallskip

We now give the intertwiner of the Schroedinger model with the Fock model for 
$\la = 2\pi i$. The $K'$-finite vectors of the Schroedinger model form the 
polynomial Fock space $S(V^n) \subset \mathcal{S}(V^n)$ which consists of those 
Schwartz functions on $V^n$ of the form $p(\x)\varphi_0(\x)$, where $p(\x)$ is 
a polynomial function on $V^n$ and $\varphi_0(\x)$ is the standard Gaussian on 
$V^n$. On the other hand, we define an action of the quantum algebra $\mathcal
{\W}_{\la}$ on $S(V^n)$ by
\begin{align*}
\omega(v_{\alpha} \otimes e_j) &= 2 \pi i x_{\alpha j},  \qquad \qquad & \omega
(v_{\alpha} \otimes f_j)  = - \frac{\partial}{\partial x_{\alpha j}}, \\ 
\omega(v_{\mu} \otimes e_j) &= - 2 \pi i x_{\mu  j},  \qquad \qquad & \omega(v_
{\mu} \otimes f_j) = - \frac{\partial}{\partial x_{\mu j}},
\end{align*}
which has central character $\la  =2\pi i$. As before, we obtain an action of 
$\mathfrak{sp}(V \otimes W)$, and this is the infinitesimal action of the 
Schroedinger model of the Weil representation introduced in the previous 
section.
For $\la = 2\pi i$, we then have a unique $\mathcal{\W}_{\la}$-intertwining 
operator $\iota: S(V^n) \rightarrow \mathcal{P}(\C^{n(p+q)})$ satisfying $\iota
(\varphi_0) = 1$ ($\W'$ annihilates $1 \in \mathcal{P}(\C^{n(p+q)})$ and 
$\varphi_0 \in S(V^n)$). 

\begin{lemma}\label{inter1}

The intertwining operator between the Schroedinger and the Fock model satisfies 
\begin{align*}
\iota
\left(  x_{\alpha j} - \frac{1}{2\pi} \frac{\partial}{\partial x_{\alpha j}}  
\right)  \iota^{-1} &= -i \frac1{2\pi} z_{\alpha j}, & 
\iota \left(  x_{\alpha j} + \frac{1}{2\pi} \frac{\partial}{\partial x_{\alpha 
j}}  \right)  \iota^{-1} &= 2  i\frac{\partial}{\partial z_{\alpha j}}, \\ 
\notag
\iota \left(  x_{\mu j} - \frac{1}{2\pi} \frac{\partial}{\partial x_{\mu j}}  
\right)  \iota^{-1} &= i \frac{1}{2\pi} z_{\mu j} , &
\iota \left(  x_{\mu j} + \frac{1}{2\pi} \frac{\partial}{\partial x_{\mu j}}  
\right)  \iota^{-1} &= -2  i\frac{\partial}{\partial z_{\mu j}}.
\end{align*}

\end{lemma}

\end{document}